\documentclass[11pt]{article}

\usepackage{amsmath}
\usepackage{amssymb}
\usepackage{amsthm}

\allowdisplaybreaks
\usepackage{parskip}
\usepackage{setspace}
\usepackage{xcolor}  
\usepackage{graphicx}
\usepackage{multicol}
\usepackage{multirow}
\usepackage[sort,round]{natbib}
\usepackage{lmodern}
\usepackage{subfig}
\usepackage{subcaption}

\newtheorem{theorem}{Theorem}[section]

\newtheorem{corollary}[theorem]{Corollary}

\newtheorem{proposition}[theorem]{Proposition}
\newtheorem{remark}[theorem]{Remark}

\usepackage{anysize}
\marginsize{15mm}{15mm}{15mm}{30mm} 

\begin{document}

\title{Asymptotic analysis of infinite-server queues driven by a Cox process with squared Bessel arrival intensity}

\author{Jos\'{e} Da Fonseca\thanks{Auckland University of Technology, Business School, Department of Economics and Finance, Private Bag 92006, 1142 Auckland, New Zealand. Phone: +64 9 9219999 extn 5063. Email: jose.dafonseca@aut.ac.nz and PRISM Sorbonne EA 4101, Universit\'e Paris 1 Panth\'eon - Sorbonne, 17 rue de la Sorbonne, 75005 Paris, France.  ORCID: 0000-0002-6882-4511} 
\and
Yajun Shi\thanks{Auckland University of Technology, School of Computing and Mathematical Sciences, Private Bag 92006, 1142 Auckland, New Zealand. Email: yajun.shi@autuni.ac.nz. ORICD: 0009-0007-5262-2724}
\and
Wenjun Zhang\thanks{Corresponding author: Auckland University of Technology, School of Computing and Mathematical Sciences, Private Bag 92006, 1142 Auckland, New Zealand. Phone: +64 9 9219999 extn 5094. Email: wenjun.zhang@aut.ac.nz. ORICD: 0000-0001-7992-8204}
}

\date{\today}
\maketitle 
\singlespacing 

\begin{abstract}
	This paper introduces a model for infinite-server queues driven by a Cox process with stochastic intensity given by a squared Bessel (BESQ) process. Exploiting the affine structure of the underlying dynamics, we derive a closed-form expression for the joint moment-generating function (MGF) of the arrival intensity, the number of customers in the system, and the integral of the number of customers in the system over a time interval. The asymptotic properties of the model are also analysed. We establish laws of large numbers, central limit theorems, large deviation principles, and moderate deviation principles for several quantities of interest. More precisely, these results are obtained for three important settings: cumulative occupancy over a finite time interval, system occupancy under a diverging long-run mean arrival intensity, and the superposition of independent infinite-server queues with Cox-BESQ arrival intensities. Numerical experiments illustrate the associated rate functions and demonstrate the effectiveness of the resulting rare-event approximations.
\end{abstract}


\singlespacing 


\section{Introduction}

The arrival process is a fundamental ingredient of stochastic queueing models, as it determines the temporal pattern by which customers, jobs, or requests enter the system. In a large part of the queueing literature, arrivals are modelled as a homogeneous Poisson process. This assumption is mathematically convenient and often leads to explicit performance formulae, limit theorems, and tractable approximations. However, it is also well known that the homogeneous Poisson assumption may be too restrictive for empirical applications. In many real-world systems, arrival rates vary over time due to seasonality, external shocks, behavioural effects, or changing environmental conditions; see, for example, \cite{HeemskerkVanLeeuwaardenMandjes2017}. Moreover, empirical studies frequently report that arrival streams display substantially more variability than is predicted by a Poisson model with deterministic intensity. In such cases, the variance of the number of arrivals over a given time interval exceeds its mean, a phenomenon commonly referred to as overdispersion.

A natural way to incorporate overdispersion is to replace the deterministic arrival rate of a Poisson process by a stochastic process. This leads to the class of Cox processes, or doubly stochastic Poisson processes; see \citet{Bremaud1981}. Conditional on the path of the stochastic intensity, arrivals follow a Poisson process, while the randomness of the intensity itself introduces an additional source of variability. By the law of total variance, the dispersion index of a Cox process, defined as the ratio of the variance to the mean of the number of arrivals, is at least one and is typically strictly greater than one. Cox processes therefore provide a flexible and analytically meaningful framework for modelling overdispersed input streams.

Examples of such overdispersed arrival mechanisms arise in many service and communication systems. Consider, for instance, emergency calls to a fire department. Under normal conditions, calls may arrive at a relatively stable rate. However, periods of high temperature, dry weather, or strong winds can substantially increase the likelihood of fires, leading to a sudden rise in call intensity. This elevated intensity may then decay gradually as conditions improve, before increasing again during another period of adverse weather. Similar effects occur in call centres, cloud-computing platforms, insurance claim systems, and social-media applications. In all these settings, the arrival intensity is not merely deterministic and time-dependent, but random, and the resulting input process can be significantly more variable than a standard Poisson process.

Queues driven by Cox processes and related stochastic-intensity arrival models have always attracted attention in the literature. In \cite{KoopsBoxmaMandjes2017}, infinite-server queues are studied under Cox input where the stochastic arrival rate is modelled as a shot-noise process with exponentially decaying jumps. Related infinite-server systems are analysed in \cite{DawPender2018} and \cite{KoopsSaxenaBoxmaMandjes2018}, where the arrival process is a Hawkes process. Hawkes processes are self-exciting point processes: arrivals increase the future arrival intensity, thereby generating clustering and overdispersion. These models are well suited to systems in which each arrival can trigger further arrivals, but they may be less appropriate when fluctuations in demand are driven by an exogenous stochastic environment rather than by endogenous self-excitation.

In this paper, we contribute to the development of queueing models by introducing an infinite-server queue driven by a Cox process whose stochastic intensity is given by a squared Bessel process (BESQ) and by analysing its tractability and asymptotic properties. The squared Bessel process is non-negative, analytically tractable, and belongs to the class of affine processes. These properties make it a natural candidate for modelling stochastic arrival intensities while preserving a high degree of mathematical tractability. The resulting model provides an alternative to Hawkes-driven infinite-server queues, such as those considered in \citet{DawPender2018} and \citet{KoopsSaxenaBoxmaMandjes2018}. In contrast to self-exciting models, the BESQ-driven Cox framework captures overdispersion through an exogenous stochastic-intensity process. At the same time, its affine structure enables explicit calculations that may be difficult to obtain in Hawkes-based frameworks.

We derive the expression for the joint moment-generating function (MGF) of the arrival intensity, the number of customers in the system and the cumulative number of customers in the system over a finite time interval for an infinite-server queue driven by a Cox-BESQ process. Further to this, we contribute to the literature by developing several asymptotic properties of queue-related variables. More precisely, we establish laws of large numbers (LLNs), central limit theorems (CLTs), large deviation principles (LDPs), and moderate deviation principles (MDPs) for several quantities of interest related to the number of customers in the system (or system occupancy). These asymptotic results are obtained for three important settings: cumulative occupancy over a finite time interval, system occupancy under a diverging long-run mean arrival intensity, and the superposition of independent infinite-server queues with Cox-BESQ arrival intensities. Numerical experiments illustrate the associated rate functions and demonstrate the effectiveness of the resulting rare-event approximations.

The remainder of the paper is organised as follows. Section~\ref{sec:model} introduces the BESQ-driven infinite-server queues, specifies the model dynamics, and derives moment formulae for the joint vector consisting of the number of customers in the system and the stochastic arrival intensity. It also reports the joint MGF of the arrival intensity, the number of customers in the system, and the cumulative number of customers in the system over a time interval. Section~\ref{sec:asymptotics} establishes three types of asymptotic results commonly analysed in the literature, namely LLNs, CLTs, LDPs, and MDPs, for three important cases: the cumulative number of customers in the system over a time interval, an infinite-server queue with a diverging arrival intensity, and sums of independent infinite-server queues. Section~\ref{sec:numerics} presents numerical examples illustrating the associated rate functions and rare-event probability approximations. Finally, Section~\ref{sec:Conclusion} concludes the paper and discusses possible directions for future research.

\section{The model}
\label{sec:model}

Let \((\Omega,\mathcal{F},\mathbb{P})\) be a complete probability space equipped with a filtration
\((\mathcal{F}_t)_{t\geq 0}\). Let
\((W_t)_{t\geq 0}\) be a standard Brownian motion and assume that the stochastic
intensity process \((\lambda_t)_{t\geq 0}\) follows the Cox--Ingersoll--Ross (CIR)
dynamics
\begin{equation}
	d\lambda_t=\kappa(\theta-\lambda_t)\,dt+\sigma\sqrt{\lambda_t}\,dW_t,
	\qquad \lambda_0\geq 0,
	\label{eq:lambda_cir}
\end{equation}
where \(\kappa>0\), \(\theta>0\), and \(\sigma>0\). The square-root process
\eqref{eq:lambda_cir}, also known as a squared Bessel-type process, is widely used in quantitative finance; see, for instance, \citet[Chap.~6]{JeanblancYorChesney2009}. The process \eqref{eq:lambda_cir} is commonly referred to in finance as the CIR process when used for interest-rate modelling and as a square-root process (due to the form of the diffusion term) in equity-derivatives modelling. However, since it is a transformed BESQ process, see \citet[Proposition~6.3.1.1]{JeanblancYorChesney2009}, we will refer to it as a BESQ process.

Conditionally on the path of \((\lambda_t)_{t\geq 0}\), let \((N_t^u)_{t\geq 0}\) be a Cox process with stochastic intensity \((\lambda_t)_{t\geq 0}\); see \citet[p.~21]{Bremaud1981}. We consider an infinite-server queue and denote by \((q_t)_{t\geq 0}\) the number of customers in the system (or number in the system). Arrivals increase \(q_t\) by one at the jump times of \(N^u\), whereas departures decrease \(q_t\) by one at the jump times of a counting process \((N_t^d)_{t\geq 0}\) with stochastic intensity \(\mu q_t\), where \(\mu>0\). Thus,
\begin{equation}
	dq_t=dN_t^u-dN_t^d.
	\label{eq:queue_dynamics}
\end{equation}
Equivalently, each customer has an independent exponential service time with
parameter \(\mu\). Throughout, we write
\[
x_t=(\lambda_t,q_t)\in \mathbb{R}_+\times\mathbb{N}.
\]

For a sufficiently regular function
\(f:\mathbb{R}_+\times\mathbb{N}\to\mathbb{R}\), the infinitesimal generator of
\((x_t)_{t\geq 0}\) is
\begin{align}
	\mathcal{L}f(\lambda,q)
	&=
	\kappa(\theta-\lambda)\partial_\lambda f(\lambda,q)
	+\frac{\sigma^2}{2}\lambda\partial_{\lambda\lambda}^2 f(\lambda,q)
	\nonumber\\
	&\quad
	+\lambda\bigl(f(\lambda,q+1)-f(\lambda,q)\bigr)
	+\mu q\bigl(f(\lambda,q-1)-f(\lambda,q)\bigr).
	\label{eq:generator}
\end{align}
Consequently, whenever the required expectations are finite, Dynkin's formula
gives
\begin{equation}
	\frac{d}{dt}\mathbb{E}\bigl[f(x_t)\bigr]
	=
	\mathbb{E}\bigl[\mathcal{L}f(x_t)\bigr].
	\label{eq:dynkin_formula}
\end{equation}

\begin{proposition}[Mixed moment equations]
	\label{prop:mixed_moment_equations}
	Let \((x_t)_{t\geq 0}=(\lambda_t,q_t)_{t\geq 0}\) have generator
	\eqref{eq:generator}. For \(p,m\in\mathbb{N}\), define
	\[
	M_{p,m}(t):=\mathbb{E}\left[\lambda_t^p q_t^m\right].
	\]
	 
	Then the mixed moment \(M_{p,m}\) satisfies
	\begin{align}
		\frac{d}{dt}M_{p,m}(t)
		&=
		\left(\kappa\theta p+\frac{\sigma^2}{2}p(p-1)\right)
		M_{p-1,m}(t)
		-\kappa p\,M_{p,m}(t)
		\nonumber\\
		&\quad
		+\sum_{i=1}^{m}\binom{m}{i}M_{p+1,m-i}(t)
		+\mu\sum_{i=2}^{m}\binom{m}{i}(-1)^i M_{p,m-i+1}(t)
		-\mu m\,M_{p,m}(t).
		\label{eq:mixed_moment_equation}
	\end{align}
\end{proposition}

\begin{proof}
	Apply \eqref{eq:dynkin_formula} to
	\(
	f(\lambda,q)=\lambda^p q^m.
	\)
	The diffusion part gives
	\[
	\kappa(\theta-\lambda)\partial_\lambda f
	+\frac{\sigma^2}{2}\lambda\partial_{\lambda\lambda}^2 f
	=
	\left(\kappa\theta p+\frac{\sigma^2}{2}p(p-1)\right)
	\lambda^{p-1}q^m
	-\kappa p\,\lambda^p q^m.
	\]
	The upward jumps yield
	\[
	\lambda\bigl(f(\lambda,q+1)-f(\lambda,q)\bigr)
	=
	\lambda^{p+1}\bigl((q+1)^m-q^m\bigr)
	=
	\sum_{i=1}^{m}\binom{m}{i}\lambda^{p+1}q^{m-i}.
	\]
	Similarly, the downward jumps give
	\begin{align*}
		\mu q\bigl(f(\lambda,q-1)-f(\lambda,q)\bigr)
		&=
		\mu\lambda^p q\bigl((q-1)^m-q^m\bigr) \\
		&=
		\mu\lambda^p q
		\left(
		\sum_{i=0}^{m}\binom{m}{i}(-1)^i q^{m-i}-q^m
		\right) \\
		&=
		-\mu m\,\lambda^p q^m
		+\mu\sum_{i=2}^{m}\binom{m}{i}(-1)^i
		\lambda^p q^{m-i+1}.
	\end{align*}
	Taking expectations gives \eqref{eq:mixed_moment_equation}.
\end{proof}

The system \eqref{eq:mixed_moment_equation} has the same triangular structure as the moment equations studied by \citet{DawPender2018} and can be analysed using
the matrix framework of \citet{DawPender2023}. It is a consequence of the affine nature of the process, see \citet{CuchieroKeller-ResselTeichmann2012}. In the present setting, however, the autonomous affine dynamics of \((\lambda_t)_{t\geq 0}\) simplify the
calculation. In particular, taking \(m=0\) in
\eqref{eq:mixed_moment_equation} gives the closed recursion
\begin{equation}
	\frac{d}{dt}\mathbb{E}[\lambda_t^k]
	=
	\left(\kappa\theta k+\frac{\sigma^2}{2}k(k-1)\right)
	\mathbb{E}[\lambda_t^{k-1}]
	-\kappa k\,\mathbb{E}[\lambda_t^k],
	\qquad k\in\mathbb{N}.
	\label{eq:lambda_moment_equation}
\end{equation}

If \((\lambda_t)_{t\geq 0}\) is stationary, then
\(\frac{d}{dt}\mathbb{E}[\lambda_t^k]=0\), and
\eqref{eq:lambda_moment_equation} yields
\begin{equation}
	\mathbb{E}[\lambda_\infty^p]
	=
	\frac{
		\prod_{i=2}^{p}
		\left(\kappa\theta+\frac{\sigma^2}{2}(i-1)\right)
	}{\kappa^{p-1}}
	\mathbb{E}[\lambda_\infty],
	\qquad
	\mathbb{E}[\lambda_\infty]=\theta.
	\label{eq:stationary_lambda_moments}
\end{equation}

\begin{proposition}[Stationary first and second moments]
	\label{prop:stationary_first_second_moments}
	Assume that \((\lambda_t,q_t)\) admits a stationary distribution.
	Then
	\begin{align}
		\mathbb{E}[q_\infty]
		&=
		\frac{\theta}{\mu},
		\label{eq:stationary_q_mean}
		\\
		\operatorname{Var}(q_\infty)
		&=
		\frac{\theta}{\mu}
		+
		\frac{\sigma^2\theta}{2\kappa\mu(\kappa+\mu)}.
		\label{eq:stationary_q_variance}
	\end{align}
	Equivalently,
	\begin{equation}
		\mathbb{E}[q_\infty^2]
		=
		\frac{\theta^2}{\mu^2}
		+
		\frac{\theta}{\mu}
		+
		\frac{\sigma^2\theta}{2\kappa\mu(\kappa+\mu)}.
		\label{eq:stationary_q_second_moment}
	\end{equation}
\end{proposition}

\begin{proof}
	Taking \(p=0,m=1\) in \eqref{eq:mixed_moment_equation} gives
	\[
	\frac{d}{dt}\mathbb{E}[q_t]
	=
	\mathbb{E}[\lambda_t]-\mu\mathbb{E}[q_t].
	\]
	At stationarity, using \(\mathbb{E}[\lambda_\infty]=\theta\), this gives
	\eqref{eq:stationary_q_mean}.
	
	Next, the equations for the variance and covariance are
	\begin{align}
		\frac{d}{dt}\operatorname{Var}(\lambda_t)
		&=
		-2\kappa\operatorname{Var}(\lambda_t)
		+\sigma^2\mathbb{E}[\lambda_t],
		\label{eq:lambda_variance_ode}
		\\
		\frac{d}{dt}\operatorname{Cov}(\lambda_t,q_t)
		&=
		-(\kappa+\mu)\operatorname{Cov}(\lambda_t,q_t)
		+\operatorname{Var}(\lambda_t),
		\label{eq:lambda_q_covariance_ode}
		\\
		\frac{d}{dt}\operatorname{Var}(q_t)
		&=
		-2\mu\operatorname{Var}(q_t)
		+2\operatorname{Cov}(\lambda_t,q_t)
		+\mathbb{E}[\lambda_t]
		+\mu\mathbb{E}[q_t].
		\label{eq:q_variance_ode}
	\end{align}
	At stationarity, \eqref{eq:lambda_variance_ode} gives
	\[
	\operatorname{Var}(\lambda_\infty)
	=
	\frac{\sigma^2\theta}{2\kappa}.
	\]
	Then \eqref{eq:lambda_q_covariance_ode} gives
	\[
	\operatorname{Cov}(\lambda_\infty,q_\infty)
	=
	\frac{\operatorname{Var}(\lambda_\infty)}{\kappa+\mu}
	=
	\frac{\sigma^2\theta}{2\kappa(\kappa+\mu)}.
	\]
	Finally, substituting these identities into \eqref{eq:q_variance_ode} yields
	\[
	\operatorname{Var}(q_\infty)
	=
	\frac{
		2\operatorname{Cov}(\lambda_\infty,q_\infty)
		+\theta+\mu\theta/\mu
	}{2\mu}
	=
	\frac{\theta}{\mu}
	+
	\frac{\sigma^2\theta}{2\kappa\mu(\kappa+\mu)}.
	\]
	Since
	\(\mathbb{E}[q_\infty^2]=\operatorname{Var}(q_\infty)+\mathbb{E}[q_\infty]^2\),
	\eqref{eq:stationary_q_second_moment} follows.
\end{proof}

The first two moments show that \(q_\infty\) is overdispersed relative to a
Poisson random variable with mean \(\theta/\mu\), except in the degenerate case
\(\sigma=0\), where the additional variance term vanishes.

\begin{corollary}[Transient first two moments and stationary autocovariance]
	\label{cor:transient_moments_autocovariance}
	The transient first two moments and stationary autocovariance of $\lambda_t$ and $q_t$ are given by
	\begin{align}
		\mathbb{E}[\lambda_t]
		&=
		\lambda_0 e^{-\kappa t}+\theta(1-e^{-\kappa t}),
		\label{eq:lambda_mean_solution}
		\\
		\mathbb{E}[q_t]
		&=
		q_0e^{-\mu t}
		+
		\int_0^t e^{-\mu(t-s)}
		\left(\lambda_0 e^{-\kappa s}+\theta(1-e^{-\kappa s})\right)\,ds \nonumber
		\\
		&=
		q_0e^{-\mu t}
		+
		(\lambda_0-\theta)\frac{e^{-\kappa t}-e^{-\mu t}}{\mu-\kappa}
		+
		\frac{\theta}{\mu}(1-e^{-\mu t}),
		\qquad \kappa\neq\mu.
		\label{eq:q_mean_solution}
	\end{align}
	Moreover,
	\begin{align}
		\operatorname{Var}(q_t)
		&=
		\int_0^t e^{-2\mu(t-s)}
		\bigl\{
		2\operatorname{Cov}(\lambda_s,q_s)
		+\mathbb{E}[\lambda_s]
		+\mu\mathbb{E}[q_s]
		\bigr\}\,ds,
		\label{eq:q_variance_integral}
	\end{align}
	where
	\begin{align}
		\operatorname{Var}(\lambda_t)
		&=
		\int_0^t
		e^{-2\kappa(t-s)}
		\sigma^2
		\left(\lambda_0e^{-\kappa s}+\theta(1-e^{-\kappa s})\right)\,ds,
		\label{eq:lambda_variance_integral}
		\\
		\operatorname{Cov}(\lambda_t,q_t)
		&=
		\int_0^t
		e^{-(\kappa+\mu)(t-s)}
		\operatorname{Var}(\lambda_s)\,ds.
		\label{eq:lambda_q_covariance_integral}
	\end{align}

	Furthermore, if \((\lambda_t,q_t)\) is stationary and \(t_2=t_1+\tau\), with
	\(\tau\geq 0\), then
	\begin{align}
		\operatorname{Cov}(q_{t_1},q_{t_2})
		&=
		e^{-\mu\tau}
		\left(
		\frac{\theta}{\mu}
		+
		\frac{\sigma^2\theta}{2\kappa\mu(\kappa+\mu)}
		\right)
		+
		\frac{\sigma^2\theta}{2\kappa(\kappa+\mu)}
		\frac{e^{-\kappa\tau}-e^{-\mu\tau}}{\mu-\kappa},
		\qquad \mu\neq\kappa.
		\label{eq:stationary_q_autocovariance}
	\end{align}
\end{corollary}

\begin{proof}
	Equations \eqref{eq:lambda_mean_solution}--\eqref{eq:q_mean_solution} follow
	from the linear system
	\begin{align}
		\frac{d}{dt}\mathbb{E}[\lambda_t]
		&=
		-\kappa\mathbb{E}[\lambda_t]+\kappa\theta,
		\label{eq:lambda_mean_ode}
		\\
		\frac{d}{dt}\mathbb{E}[q_t]
		&=
		-\mu\mathbb{E}[q_t]+\mathbb{E}[\lambda_t],
		\label{eq:q_mean_ode}
	\end{align}
	obtained from \eqref{eq:mixed_moment_equation}. Similarly,
	\eqref{eq:lambda_variance_integral} and \eqref{eq:lambda_q_covariance_integral} follow by
	variation of constants applied to
	\eqref{eq:lambda_variance_ode}--\eqref{eq:q_variance_ode}. 
	
	It remains to identify the stationary autocovariance. Conditional on
	\(\mathcal{F}_{t_1}\), the conditional mean of \(q_{t_2}\) is
	\[
	\mathbb{E}[q_{t_2}\mid\mathcal{F}_{t_1}]
	=
	e^{-\mu\tau}q_{t_1}
	+
	\int_0^\tau e^{-\mu(\tau-s)}
	\mathbb{E}[\lambda_{t_1+s}\mid\mathcal{F}_{t_1}]\,ds.
	\]
	Since the CIR process is Markovian,
	\[
	\mathbb{E}[\lambda_{t_1+s}\mid\mathcal{F}_{t_1}]
	=
	\theta+(\lambda_{t_1}-\theta)e^{-\kappa s}.
	\]
	Therefore,
	\begin{align*}
		\operatorname{Cov}(q_{t_1},q_{t_2})
		&=
		e^{-\mu\tau}\operatorname{Var}(q_{t_1})
		+
		\operatorname{Cov}(q_{t_1},\lambda_{t_1})
		\int_0^\tau e^{-\mu(\tau-s)}e^{-\kappa s}\,ds.
	\end{align*}
	where under stationarity, the explicit forms of the terms are specified as 
	\[
	\operatorname{Var}(q_{t_1})
	=
	\frac{\theta}{\mu}
	+
	\frac{\sigma^2\theta}{2\kappa\mu(\kappa+\mu)},
	\qquad
	\operatorname{Cov}(q_{t_1},\lambda_{t_1})
	=
	\frac{\sigma^2\theta}{2\kappa(\kappa+\mu)}.
	\]
	and
	\[
	\int_0^\tau e^{-\mu(\tau-s)}e^{-\kappa s}\,ds
	=
	\frac{e^{-\kappa\tau}-e^{-\mu\tau}}{\mu-\kappa},
	\qquad \mu\neq\kappa.
	\]
\end{proof}


A contribution of this work is the derivation of the MGF of the incoming intensity, the number of customers in the system, and the cumulative number of customers in the system over a given time interval. Verifying that this function is well defined is achieved by comparing the ODEs involved in the MGF with the solution of a scalar Riccati equation with constant coefficients. These two results, \textit{i.e.}, a comparison property for solutions of Riccati ODEs and the well-definedness of a scalar Riccati equation, are both well known and are provided for completeness.

\begin{proposition}[Explicit solution and explosion of a constant-coefficient Riccati equation]
	\label{prop:riccati_constant_solution_explosion}
	Let \((y(t))_{t\geq 0}\) solve
	\begin{align}
		y'(t)
		&=
		\frac{\sigma^2}{2}y(t)^2-\kappa y(t)+c,
		\label{eq:riccati_constant_ode}
		\\
		y(0)
		&=
		u,
		\label{eq:riccati_constant_initial_condition}
	\end{align}
	where \(\sigma>0\), \(\kappa>0\), \(u\in\mathbb{R}\), and
	\(c\in\mathbb{R}\). Assume that
	\[
	c\leq \frac{\kappa^2}{2\sigma^2}.
	\]
	Define
	\begin{equation}
		\Delta:=\kappa^2-2\sigma^2c\geq 0,
		\qquad
		\nu_\pm:=\frac{-\kappa\pm\sqrt{\Delta}}{2},
		\qquad
		\nu_-<\nu_+.
		\label{eq:riccati_delta_nu}
	\end{equation}
	If \(\Delta>0\), then, up to its explosion time,
	\begin{equation}
		y(t)
		=
		-\frac{
			(c+u\nu_+)-(c+u\nu_-)e^{-\sqrt{\Delta}t}
		}{
			\left(\nu_-+\frac{\sigma^2}{2}u\right)
			-
			\left(\nu_++\frac{\sigma^2}{2}u\right)e^{-\sqrt{\Delta}t}
		}.
		\label{eq:riccati_constant_solution}
	\end{equation}
	The explosion time
	\[
	t^\star:=\inf\{t>0:y(t)=+\infty\},
	\]
	when it exists, is
	\begin{equation}
		t^\star
		=
		\frac{1}{\sqrt{\Delta}}
		\log\left(
		\frac{2\nu_+ +\sigma^2 u}{2\nu_-+\sigma^2 u}
		\right).
		\label{eq:riccati_constant_explosion_time}
	\end{equation}
	In particular, there is no positive explosion time whenever
	\begin{equation}
		u<\frac{\kappa+\sqrt{\Delta}}{\sigma^2}.
		\label{eq:riccati_constant_no_explosion_condition}
	\end{equation}
\end{proposition}

\begin{proof}
	Use the standard Riccati transformation
	\begin{equation}
		y(t)=-\frac{h'(t)}{\frac{\sigma^2}{2}h(t)}.
		\label{eq:riccati_linearisation}
	\end{equation}
	Substituting \eqref{eq:riccati_linearisation} into
	\eqref{eq:riccati_constant_ode} gives the second-order linear ODE
	\begin{equation}
		h''(t)+\kappa h'(t)+\frac{\sigma^2 c}{2}h(t)=0.
		\label{eq:riccati_linear_ode}
	\end{equation}
	The characteristic polynomial is
	\[
	\nu^2+\kappa\nu+\frac{\sigma^2c}{2}=0,
	\]
	whose roots are \(\nu_-\) and \(\nu_+\), as defined in
	\eqref{eq:riccati_delta_nu}. Hence, when \(\Delta>0\),
	\[
	h(t)=C_-e^{\nu_-t}+C_+e^{\nu_+t}.
	\]
	The initial condition \(y(0)=u\), together with
	\eqref{eq:riccati_linearisation}, implies
	\[
	\frac{C_-}{C_+}
	=
	-\frac{2\nu_+ +\sigma^2u}{2\nu_-+\sigma^2u}.
	\]
	Substituting this ratio into \eqref{eq:riccati_linearisation} yields
	\eqref{eq:riccati_constant_solution}.
	
	Explosion occurs precisely when the denominator in
	\eqref{eq:riccati_constant_solution} vanishes. Thus \(t^\star\) must satisfy
	\[
	\left(\nu_-+\frac{\sigma^2}{2}u\right)
	=
	\left(\nu_++\frac{\sigma^2}{2}u\right)e^{-\sqrt{\Delta}t^\star},
	\]
	which gives \eqref{eq:riccati_constant_explosion_time}, provided that the
	right-hand side defines a strictly positive time. Since
	\[
	2\nu_+ + \sigma^2u
	=
	-\kappa+\sqrt{\Delta}+\sigma^2u,
	\qquad
	2\nu_- + \sigma^2u
	=
	-\kappa-\sqrt{\Delta}+\sigma^2u,
	\]
	a positive explosion time can only occur when
	\[
	u>\frac{\kappa+\sqrt{\Delta}}{\sigma^2}.
	\]
	Therefore no positive explosion time exists whenever
	\eqref{eq:riccati_constant_no_explosion_condition} holds.
\end{proof}

\begin{corollary}[Comparison for Riccati equations]
	\label{cor:riccati_comparison}
	Let \((y(t))_{t\geq 0}\) solve \eqref{eq:riccati_constant_ode} with constant
	source term \(c\), and let \((b(t))_{t\geq 0}\) solve
	\begin{align}
		b'(t)
		&=
		\frac{\sigma^2}{2}b(t)^2-\kappa b(t)+c(t),
		\label{eq:riccati_time_dependent_ode}
		\\
		b(0)
		&=
		u,
		\label{eq:riccati_time_dependent_initial_condition}
	\end{align}
	where \(c(t)\leq c\) for all \(t\geq 0\). Then, up to the first explosion time
	of either equation,
	\begin{equation}
		y(t)\geq b(t).
		\label{eq:riccati_comparison_result}
	\end{equation}
\end{corollary}

\begin{proof}
	Let
	\[
	z(t):=y(t)-b(t).
	\]
	Then \(z(0)=0\) and
	\begin{align*}
		z'(t)
		&=
		y'(t)-b'(t) \\
		&=
		\left(
		\frac{\sigma^2}{2}(y(t)+b(t))-\kappa
		\right)z(t)
		+
		c-c(t).
	\end{align*}
	Set
	\[
	A(t):=\frac{\sigma^2}{2}(y(t)+b(t))-\kappa.
	\]
	By variation of constants,
	\[
	z(t)
	=
	e^{\int_0^t A(s)\,ds}
	\left(
	z(0)
	+
	\int_0^t
	(c-c(r))e^{-\int_0^r A(s)\,ds}\,dr
	\right).
	\]
	Since \(z(0)=0\) and \(c-c(r)\geq 0\), it follows that \(z(t)\geq 0\). Hence
	\(y(t)\geq b(t)\).
\end{proof}


A key property of the infinite-server queue with Cox--BESQ arrivals is that the
joint moment-generating function of
\[
\left(\lambda_T,q_T,\int_t^T q_s\,ds\right)
\]
has an explicit exponentially affine form. This follows from the affine structure of the Markov process \((\lambda_t,q_t)_{t\geq 0}\) and from the explicit solvability of the Riccati equations involved. Such explicit solvability is a useful feature of the present model and is not available in some affine point-process models, such as general Hawkes-type specifications (for a detailed example, see \citet{DaFonsecaZaatour2014}).

\begin{theorem}[Joint MGF of \((\lambda_T, q_T,\int_t^T q_s ds)\)]
	\label{thm:mgf_triplet}
	Let \(0\leq t<T\) and let
	\(\mathbf{u}=(u_1,u_2,u_3)\in\mathbb{R}^3\). Define
	\begin{equation}
		f(t,T,\lambda_t,q_t)
		:=
		\mathbb{E}_t\left[
		\exp\left(
		u_1\lambda_T
		+u_2q_T
		+u_3\int_t^T q_s\,ds
		\right)
		\right].
		\label{eq:mgf_triplet_expectation}
	\end{equation}
	Then, on the parameter domain where the expectation is finite,
	\begin{equation}
		f(t,T,\lambda_t,q_t)
		=
		\exp\left(
		a(T-t;\mathbf{u})
		+b(T-t;\mathbf{u})\lambda_t
		+d(T-t;\mathbf{u})q_t
		\right).
		\label{eq:mgf_triplet_affine_form}
	\end{equation}
	The function \(d\) is given by
	\begin{equation}
		d(s;\mathbf{u})
		=
		\log\left(
		e^{(u_3-\mu)s}e^{u_2}
		-
		\frac{\mu}{u_3-\mu}
		\left(1-e^{(u_3-\mu)s}\right)
		\right),
		\qquad s\geq 0.
		\label{eq:mgf_triplet_d_solution}
	\end{equation}

	Let
	\begin{equation}
		\alpha
		:=
		\frac{\sigma^2}{2}\frac{u_3}{u_3-\mu},
		\qquad
		\beta
		:=
		\frac{\sigma^2}{2}
		\left(
		e^{u_2}+\frac{\mu}{u_3-\mu}
		\right),
		\qquad
		k_1
		:=
		\frac{\sqrt{\kappa^2+4\alpha}}{\mu-u_3}.
		\label{eq:mgf_triplet_parameters}
	\end{equation}
	Assume that \(\kappa^2+4\alpha>0\). Define
	\begin{align}
		h_1(s)
		&:=
		e^{\frac{-\kappa-\sqrt{\kappa^2+4\alpha}}{2}s}
		{}_0F_1\left(
		1+k_1,
		-\frac{\beta e^{(u_3-\mu)s}}{(\mu-u_3)^2}
		\right),
		\label{eq:mgf_triplet_h1}
		\\
		h_2(s)
		&:=
		e^{\frac{-\kappa+\sqrt{\kappa^2+4\alpha}}{2}s}
		{}_0F_1\left(
		1-k_1,
		-\frac{\beta e^{(u_3-\mu)s}}{(\mu-u_3)^2}
		\right),
		\label{eq:mgf_triplet_h2}
	\end{align}
	where \({}_0F_1\) denotes the generalized hypergeometric function; see
	\citet[(16.2.1)]{DLMF}. Their derivatives are
	\begin{align}
		h_1'(s)
		&=
		\frac{-\kappa-\sqrt{\kappa^2+4\alpha}}{2}h_1(s)
		+
		\frac{
			e^{\frac{-\kappa-\sqrt{\kappa^2+4\alpha}}{2}s}
		}{1+k_1}
		\frac{\beta e^{(u_3-\mu)s}}{|\mu-u_3|}
		{}_0F_1\left(
		2+k_1,
		-\frac{\beta e^{(u_3-\mu)s}}{(\mu-u_3)^2}
		\right),
		\label{eq:mgf_triplet_h1_prime}
		\\
		h_2'(s)
		&=
		\frac{-\kappa+\sqrt{\kappa^2+4\alpha}}{2}h_2(s)
		+
		\frac{
			e^{\frac{-\kappa+\sqrt{\kappa^2+4\alpha}}{2}s}
		}{1-k_1}
		\frac{\beta e^{(u_3-\mu)s}}{|\mu-u_3|}
		{}_0F_1\left(
		2-k_1,
		-\frac{\beta e^{(u_3-\mu)s}}{(\mu-u_3)^2}
		\right).
		\label{eq:mgf_triplet_h2_prime}
	\end{align}
	Then
	\begin{align}
		b(s;\mathbf{u})
		&=
		-\frac{2}{\sigma^2}
		\frac{
			\left(\frac{2}{\sigma^2}h_1'(0)+u_1h_1(0)\right)h_2'(s)
			-
			\left(\frac{2}{\sigma^2}h_2'(0)+u_1h_2(0)\right)h_1'(s)
		}{
			\left(\frac{2}{\sigma^2}h_1'(0)+u_1h_1(0)\right)h_2(s)
			-
			\left(\frac{2}{\sigma^2}h_2'(0)+u_1h_2(0)\right)h_1(s)
		},
		\label{eq:mgf_triplet_b_solution}
		\\
		a(s;\mathbf{u})
		&=
		-\frac{2\kappa\theta}{\sigma^2}
		\log\left(
		\frac{
			\left(\frac{2}{\sigma^2}h_1'(0)+u_1h_1(0)\right)h_2(s)
			-
			\left(\frac{2}{\sigma^2}h_2'(0)+u_1h_2(0)\right)h_1(s)
		}{
			\left(\frac{2}{\sigma^2}h_1'(0)+u_1h_1(0)\right)h_2(0)
			-
			\left(\frac{2}{\sigma^2}h_2'(0)+u_1h_2(0)\right)h_1(0)
		}
		\right).
		\label{eq:mgf_triplet_a_solution}
	\end{align}
	Furthermore, there exist \(t^\star>0\) and positive thresholds
	\(u_1^\star,u_2^\star,u_3^\star>0\) such that
	\eqref{eq:mgf_triplet_affine_form} is well defined whenever
	\(u_i<u_i^\star\), \(i=1,2,3\), and \(T-t<t^\star\).
\end{theorem}

\begin{proof}
	By the Feynman--Kac formula, \(f\) satisfies
	\begin{equation}
		\partial_t f(t,T,\lambda,q)
		+
		\mathcal{L}f(t,T,\lambda,q)
		+
		u_3 q f(t,T,\lambda,q)
		=
		0,
		\qquad
		f(T,T,\lambda,q)=e^{u_1\lambda+u_2q}.
		\label{eq:mgf_triplet_pde}
	\end{equation}
	We use the exponentially affine ansatz
	\[
	f(t,T,\lambda,q)
	=
	\exp\left(a(s)+b(s)\lambda+d(s)q\right),
	\qquad
	s:=T-t,
	\]
	which corresponds to \eqref{eq:mgf_triplet_affine_form}, with the dependence on \(\mathbf{u}\) omitted. 	Substitution into \eqref{eq:mgf_triplet_pde}, followed by identification of
	the constant, \(\lambda\), and \(q\) coefficients, gives
	\begin{align}
		a'(s)
		&=
		\kappa\theta b(s),
		&
		a(0)&=0,
		\label{eq:mgf_triplet_a_ode}
		\\
		b'(s)
		&=
		\frac{\sigma^2}{2}b(s)^2-\kappa b(s)+e^{d(s)}-1,
		&
		b(0)&=u_1,
		\label{eq:mgf_triplet_b_ode}
		\\
		d'(s)
		&=
		\mu(e^{-d(s)}-1)+u_3,
		&
		d(0)&=u_2.
		\label{eq:mgf_triplet_d_ode}
	\end{align}
	
	Solving \eqref{eq:mgf_triplet_d_ode} gives
	\eqref{eq:mgf_triplet_d_solution}. Therefore
	\[
	e^{d(s)}-1	= \frac{2}{\sigma^2}(\beta e^{(u_3-\mu)s}-\alpha),
	\]
	with \(\alpha\) and \(\beta\) as in \eqref{eq:mgf_triplet_parameters}. Hence
	\eqref{eq:mgf_triplet_b_ode} becomes
	\begin{equation}
		b'(s)
		=
		\frac{\sigma^2}{2}b(s)^2-\kappa b(s)
		+\frac{2}{\sigma^2}(\beta e^{(u_3-\mu)s}-\alpha).
		\label{eq:mgf_triplet_b_ode_reduced}
	\end{equation}
	Using the Riccati transformation
	\begin{align}
	b(s)=-\frac{h'(s)}{\frac{\sigma^2}{2}h(s)}, \label{eq:Riccati_linearized}
	\end{align}
	we obtain the second-order linear ODE
	\begin{equation}
		h''(s)+\kappa h'(s)
		+\left(\beta e^{(u_3-\mu)s}	- \alpha\right)h(s)=0.
		\label{eq:mgf_triplet_h_ode}
	\end{equation}
	This equation has the two independent solutions \eqref{eq:mgf_triplet_h1} and
	\eqref{eq:mgf_triplet_h2}. These are obtained by setting
	\[
	h(s)=e^{c_0s}g\left(-\frac{\beta e^{(u_3-\mu)s}}{(\mu-u_3)^2}\right),
	\]
	and choosing
	\[
	c_0=\frac{-\kappa\pm\sqrt{\kappa^2+4\alpha}}{2}
	\]
	leads for the equation for \(g\) to
	\[
	zg''(z)+(1\pm k_1)g'(z)-g(z)=0,
	\]
	whose solution is expressed in terms of \({}_0F_1\). The derivative formulas
	\eqref{eq:mgf_triplet_h1_prime} and \eqref{eq:mgf_triplet_h2_prime} follow from
	the identity
	\[
	\frac{d}{dz}{}_0F_1(\nu,z)=\frac{1}{\nu}{}_0F_1(\nu+1,z).
	\]
	
	Thus the general solution for \(b\) is
	\[
	b(s)
	=
	-\frac{2}{\sigma^2}
	\frac{\bar C_1h_1'(s)+\bar C_2h_2'(s)}
	{\bar C_1h_1(s)+\bar C_2h_2(s)}.
	\]
	The initial condition \(b(0)=u_1\) implies
	\[
	\frac{\bar C_1}{\bar C_2}
	=
	-
	\frac{\frac{2}{\sigma^2}h_2'(0)+u_1h_2(0)}
	{\frac{2}{\sigma^2}h_1'(0)+u_1h_1(0)}.
	\]
	Substitution gives \eqref{eq:mgf_triplet_b_solution}. Finally, since
	\(a'(s)=\kappa\theta b(s)\) and \(b\) is a logarithmic derivative, integration
	gives \eqref{eq:mgf_triplet_a_solution}.
	
	It remains to justify that the expression is finite on a non-empty parameter
	domain. The condition \(\kappa^2+4\alpha>0\) is satisfied, for instance, when
	\[
	u_3<\frac{\mu}{1+2\sigma^2/\kappa^2}.
	\]
	Moreover, writing \(d(s)=\log g(u_2,u_3,s)\), where
	\[
	g(u_2,u_3,s)
	=
	e^{(u_3-\mu)s}e^{u_2}
	-
	\frac{\mu}{u_3-\mu}\left(1-e^{(u_3-\mu)s}\right),
	\]
	we see that \(g>0\) for all sufficiently small \(u_2,u_3\) and \(s\). Hence
	\(d\) is well defined locally. Since \(e^{d(s)}-1\) is locally bounded, the
	comparison result in Corollary~\ref{cor:riccati_comparison}, together with
	Proposition~\ref{prop:riccati_constant_solution_explosion}, implies that
	\(b\) is finite on a sufficiently small time interval. The function \(a\) is
	then finite because it is obtained by integrating \(b\). This proves the stated
	local well-posedness of the MGF.
\end{proof}

The equation \eqref{eq:Riccati_linearized} shows that the Riccati equation can be linearized and is therefore solvable (see, e.g., \citet{GrasselliTebaldi2008}). One remarkable property of this Riccati ODE is that, despite one of its coefficients being time-dependent, it can be integrated explicitly.

The MGF of \(\lambda_T\) follows from Theorem~\ref{thm:mgf_triplet} by setting \(\mathbf{u}=(u,0,0)\). This classical formula can also be found, for example,
in \citet[Cor.~6.3.4.4]{JeanblancYorChesney2009}.

\begin{corollary}[MGF of the CIR intensity]
	\label{cor:mgf_lambda}
	For \(0\leq t<T\) and \(u\in\mathbb{R}\), define
	\begin{equation}
		f(t,T,\lambda_t)
		:=
		\mathbb{E}_t\left[e^{u\lambda_T}\right].
		\label{eq:mgf_lambda_expectation}
	\end{equation}
	Then, on its domain of finiteness,
	\begin{equation}
		f(t,T,\lambda_t)
		=
		\exp\left(a(T-t;u)+b(T-t;u)\lambda_t\right),
		\label{eq:mgf_lambda_affine_form}
	\end{equation}
	where
	\begin{align}
		b(s;u)
		&=
		\frac{ue^{-\kappa s}}
		{1-\frac{u\sigma^2}{2\kappa}(1-e^{-\kappa s})},
		\label{eq:mgf_lambda_b_solution}
		\\
		a(s;u)
		&=
		-\frac{2\kappa\theta}{\sigma^2}
		\log\left(
		1-\frac{u\sigma^2}{2\kappa}(1-e^{-\kappa s})
		\right).
		\label{eq:mgf_lambda_a_solution}
	\end{align}
	For \(u>2\kappa/\sigma^2\), the explosion time is
	\begin{equation}
		t^\star(u)
		=
		-\frac{1}{\kappa}
		\log\left(
		1-\frac{2\kappa}{u\sigma^2}
		\right).
		\label{eq:mgf_lambda_explosion_time}
	\end{equation}
\end{corollary}

\begin{proof}
	Set \(\mathbf{u}=(u,0,0)\) in Theorem~\ref{thm:mgf_triplet} and, wherever it appears, replace \(\mathbf{u}\) with \(u\). Then
	\(d(s)=0\), and the Riccati system reduces to
	\[
	b'(s)=\frac{\sigma^2}{2}b(s)^2-\kappa b(s),
	\qquad
	b(0)=u,
	\]
	with \(a'(s)=\kappa\theta b(s)\) and \(a(0)=0\). Solving this scalar Riccati
	equation gives \eqref{eq:mgf_lambda_b_solution}, and integration gives
	\eqref{eq:mgf_lambda_a_solution}. Explosion occurs when the denominator in
	\eqref{eq:mgf_lambda_b_solution} vanishes, which gives
	\eqref{eq:mgf_lambda_explosion_time}.
\end{proof}

\begin{corollary}[MGF of \(q_T\)]
	\label{cor:mgf_q}
	Assume that \(\kappa/\mu\notin\mathbb{N}\). For \(0\leq t<T\) and
	\(u\in\mathbb{R}\), define
	\begin{equation}
		f(t,T,\lambda_t,q_t)
		:=
		\mathbb{E}_t\left[e^{uq_T}\right].
		\label{eq:mgf_q_expectation}
	\end{equation}
	Then, on its domain of finiteness,
	\begin{equation}
		f(t,T,\lambda_t,q_t)
		=
		\exp\left(
		a(T-t;u)+b(T-t;u)\lambda_t+d(T-t;u)q_t
		\right),
		\label{eq:mgf_q_affine_form}
	\end{equation}
	where
	\begin{align}
		d(s;u)
		&=
		\log\left(1+e^{-\mu s}(e^u-1)\right),
		\label{eq:mgf_q_d_solution}
		\\
		b(s;u)
		&=
		-\frac{2}{\sigma^2}
		\frac{h_1'(0)h_2'(s)-h_2'(0)h_1'(s)}
		{h_1'(0)h_2(s)-h_2'(0)h_1(s)},
		\label{eq:mgf_q_b_solution}
		\\
		a(s;u)
		&=
		-\frac{2\kappa\theta}{\sigma^2}
		\log\left(
		\frac{h_1'(0)h_2(s)-h_2'(0)h_1(s)}
		{h_1'(0)h_2(0)-h_2'(0)h_1(0)}
		\right).
		\label{eq:mgf_q_a_solution}
	\end{align}
	Here
	\begin{equation}
		\beta:=\frac{\sigma^2}{2}(e^u-1),
		\label{eq:mgf_q_beta}
	\end{equation}
	and
	\begin{align}
		h_1(s)
		&=
		e^{-\kappa s}
		{}_0F_1\left(
		1+\frac{\kappa}{\mu},
		-\frac{\beta}{\mu^2}e^{-\mu s}
		\right),
		\label{eq:mgf_q_h1}
		\\
		h_2(s)
		&=
		{}_0F_1\left(
		1-\frac{\kappa}{\mu},
		-\frac{\beta}{\mu^2}e^{-\mu s}
		\right),
		\label{eq:mgf_q_h2}
		\\
		h_1'(s)
		&=
		-\kappa h_1(s)
		+
		\frac{e^{-(\kappa+\mu)s}}{1+\kappa/\mu}
		\frac{\beta}{\mu}
		{}_0F_1\left(
		2+\frac{\kappa}{\mu},
		-\frac{\beta}{\mu^2}e^{-\mu s}
		\right),
		\label{eq:mgf_q_h1_prime}
		\\
		h_2'(s)
		&=
		\frac{e^{-\mu s}}{1-\kappa/\mu}
		\frac{\beta}{\mu}
		{}_0F_1\left(
		2-\frac{\kappa}{\mu},
		-\frac{\beta}{\mu^2}e^{-\mu s}
		\right).
		\label{eq:mgf_q_h2_prime}
	\end{align}
	Moreover, there exist \(t^\star>0\) and \(u^\star>0\) such that
	\eqref{eq:mgf_q_affine_form} is well defined whenever \(u<u^\star\) and
	\(T-t<t^\star\).
\end{corollary}

\begin{proof}
	Set \(\mathbf{u}=(0,u,0)\) in Theorem~\ref{thm:mgf_triplet} and, wherever it appears, replace \(\mathbf{u}\) with \(u\). Then
	\[
	\alpha=0,
	\qquad
	k_1=\frac{\kappa}{\mu},
	\qquad
	\beta=\frac{\sigma^2}{2}(e^u-1),
	\]
	and \eqref{eq:mgf_q_d_solution}--\eqref{eq:mgf_q_h2_prime} follow directly
	from the formulas in Theorem~\ref{thm:mgf_triplet}.
	
	For the finiteness statement, write
	\[
	d(s;u)=\log g(u,s),
	\qquad
	g(u,s)=1+e^{-\mu s}(e^u-1).
	\]
	For all \(u\in\mathbb{R}\) and \(s\geq 0\), \(g(u,s)>0\), so \(d\) is well
	defined. Since \(u\mapsto g(u,s)-1\) is increasing, there exist
	\(u^\star>0\) and \(t^\star>0\) such that
	\[
	g(u,s)-1<\frac{\kappa^2}{2\sigma^2},
	\qquad
	u<u^\star,
	\quad
	s<t^\star.
	\]
	The Riccati equation for \(b\) can then be compared with the constant-source
	Riccati equation of Proposition~\ref{prop:riccati_constant_solution_explosion}.
	By Corollary~\ref{cor:riccati_comparison}, \(b\) does not explode on
	\([0,t^\star)\). Since \(a(s)=\kappa\theta\int_0^s b(r)\,dr\), \(a\) is also
	finite on this interval.
\end{proof}

\begin{corollary}[Long term MGF of \(q_T\)]
	\label{cor:mgf_q_stationary}
	Assume that \(\kappa/\mu\notin\mathbb{N}\), and define
	\begin{equation}
		f_\infty(u)
		:=
		\lim_{T\to\infty}
		\mathbb{E}_t\left[e^{uq_T}\right],
		\label{eq:mgf_q_stationary_definition}
	\end{equation}
	whenever the limit exists. Then, on its domain of finiteness,
	\begin{equation}
		f_\infty(u)=e^{a(\infty;u)},
		\label{eq:mgf_q_stationary_solution}
	\end{equation}
	where
	\begin{equation}
		a(\infty;u)
		=
		\frac{2\kappa\theta}{\sigma^2}
		\log\left(
		h_2(0)
		-
		\frac{h_2'(0)h_1(0)}{h_1'(0)}
		\right),
		\label{eq:mgf_q_stationary_a_infty}
	\end{equation}
	with \(h_1,h_2,h_1'\), and \(h_2'\) given in
	Corollary~\ref{cor:mgf_q}. Moreover, as \(\beta\to 0\),
	\begin{equation}
		a(\infty;u)
		=
		\frac{\theta}{\mu}(e^u-1)+O(\beta^2),
		\qquad
		\beta=\frac{\sigma^2}{2}(e^u-1).
		\label{eq:mgf_q_stationary_poisson_expansion}
	\end{equation}
\end{corollary}

\begin{proof}
	From \eqref{eq:mgf_q_d_solution},
	\[
	\lim_{s\to\infty}d(s;u)=0.
	\]
	Furthermore, since \(\kappa>0\), \(\mu>0\), and
	\({}_0F_1(a,0)=1\),
	\[
	\lim_{s\to\infty}h_1(s)=0,
	\qquad
	\lim_{s\to\infty}h_1'(s)=0,
	\qquad
	\lim_{s\to\infty}h_2(s)=1,
	\qquad
	\lim_{s\to\infty}h_2'(s)=0.
	\]
	It follows from \eqref{eq:mgf_q_b_solution} that
	\[
	\lim_{s\to\infty}b(s;u)=0.
	\]
	Taking the limit \(s\to\infty\) in \eqref{eq:mgf_q_a_solution} gives
	\eqref{eq:mgf_q_stationary_a_infty}. The expansion
	\eqref{eq:mgf_q_stationary_poisson_expansion} follows by expanding the
	hypergeometric functions in powers of \(\beta\).
\end{proof}

\begin{remark}
	\label{rem:stationary_q_poisson_approximation}
	Equation \eqref{eq:mgf_q_stationary_poisson_expansion} shows that, to first
	order in \(\beta\), the stationary queue length \(q_\infty\) has the same MGF
	as a Poisson random variable with parameter \(\theta/\mu\). The higher-order
	terms encode the overdispersion induced by the stochastic intensity.
\end{remark}

\begin{corollary}[MGF of \(\int_t^T q_s ds\)]
	\label{cor:mgf_integrated_q}
	Let \(0\leq t<T\), \(u\in\mathbb{R}\), and assume
	\[
	u<\frac{\mu}{1+2\sigma^2/\kappa^2}.
	\]
	Define
	\begin{equation}
		f(t,T,\lambda_t,q_t)
		:=
		\mathbb{E}_t\left[
		\exp\left(u\int_t^T q_s\,ds\right)
		\right].
		\label{eq:mgf_integrated_q_expectation}
	\end{equation}
	Then, on its domain of finiteness,
	\begin{equation}
		f(t,T,\lambda_t,q_t)
		=
		\exp\left(
		a(T-t;u)+b(T-t;u)\lambda_t+d(T-t;u)q_t
		\right),
		\label{eq:mgf_integrated_q_affine_form}
	\end{equation}
	where
	\begin{equation}
		d(s;u)
		=
		\log\left(
		1-\frac{u}{u-\mu}\left(1-e^{(u-\mu)s}\right)
		\right),
		\label{eq:mgf_integrated_q_d_solution}
	\end{equation}
	and
	\begin{align}
		b(s;u)
		&=
		-\frac{2}{\sigma^2}
		\frac{h_1'(0)h_2'(s)-h_2'(0)h_1'(s)}
		{h_1'(0)h_2(s)-h_2'(0)h_1(s)},
		\label{eq:mgf_integrated_q_b_solution}
		\\
		a(s;u)
		&=
		-\frac{2\kappa\theta}{\sigma^2}
		\log\left(
		\frac{h_1'(0)h_2(s)-h_2'(0)h_1(s)}
		{h_1'(0)h_2(0)-h_2'(0)h_1(0)}
		\right).
		\label{eq:mgf_integrated_q_a_solution}
	\end{align}
	Here
	\begin{equation}
		\alpha
		=
		\frac{\sigma^2}{2}\frac{u}{u-\mu},
		\qquad
		k_1
		=
		\frac{\sqrt{\kappa^2+4\alpha}}{\mu-u},
		\label{eq:mgf_integrated_q_parameters}
	\end{equation}
	and
	\begin{align}
		h_1(s)
		&=
		e^{\frac{-\kappa-\sqrt{\kappa^2+4\alpha}}{2}s}
		{}_0F_1\left(
		1+k_1,
		-\frac{\alpha e^{(u-\mu)s}}{(\mu-u)^2}
		\right),
		\label{eq:mgf_integrated_q_h1}
		\\
		h_2(s)
		&=
		e^{\frac{-\kappa+\sqrt{\kappa^2+4\alpha}}{2}s}
		{}_0F_1\left(
		1-k_1,
		-\frac{\alpha e^{(u-\mu)s}}{(\mu-u)^2}
		\right),
		\label{eq:mgf_integrated_q_h2}
		\\
		h_1'(s)
		&=
		\frac{-\kappa-\sqrt{\kappa^2+4\alpha}}{2}h_1(s)
		+
		\frac{
			e^{\frac{-\kappa-\sqrt{\kappa^2+4\alpha}}{2}s}
		}{1+k_1}
		\frac{\alpha e^{(u-\mu)s}}{|\mu-u|}
		{}_0F_1\left(
		2+k_1,
		-\frac{\alpha e^{(u-\mu)s}}{(\mu-u)^2}
		\right),
		\label{eq:mgf_integrated_q_h1_prime}
		\\
		h_2'(s)
		&=
		\frac{-\kappa+\sqrt{\kappa^2+4\alpha}}{2}h_2(s)
		+
		\frac{
			e^{\frac{-\kappa+\sqrt{\kappa^2+4\alpha}}{2}s}
		}{1-k_1}
		\frac{\alpha e^{(u-\mu)s}}{|\mu-u|}
		{}_0F_1\left(
		2-k_1,
		-\frac{\alpha e^{(u-\mu)s}}{(\mu-u)^2}
		\right).
		\label{eq:mgf_integrated_q_h2_prime}
	\end{align}
	Moreover, there exist \(t^\star>0\) and \(u^\star>0\) such that
	\eqref{eq:mgf_integrated_q_affine_form} is well defined whenever
	\(u<u^\star\) and \(T-t<t^\star\).
\end{corollary}

\begin{proof}
	Set \(\mathbf{u}=(0,0,u)\) in Theorem~\ref{thm:mgf_triplet} and, wherever it appears, replace \(\mathbf{u}\) with \(u\). Then
	\[
	\alpha=\beta=\frac{\sigma^2}{2}\frac{u}{u-\mu},
	\]
	and the formulas follow directly from
	\eqref{eq:mgf_triplet_d_solution}--\eqref{eq:mgf_triplet_a_solution}.
	
	For the local finiteness statement, let
	\[
	g(u,s)
	:=
	1-\frac{u}{u-\mu}\left(1-e^{(u-\mu)s}\right).
	\]
	For \(u<\mu\), \(g(u,s)>0\) for all \(s\geq 0\), so \(d(s;u)=\log g(u,s)\) is
	well defined. Moreover, for sufficiently small \(u\) and \(s\), the source term
	\(e^{d(s;u)}-1=g(u,s)-1\) is bounded above by
	\(\kappa^2/(2\sigma^2)\). Corollary~\ref{cor:riccati_comparison} and
	Proposition~\ref{prop:riccati_constant_solution_explosion} therefore imply
	that \(b\) is finite on a sufficiently small time interval. Since
	\(a(s)=\kappa\theta\int_0^s b(r)\,dr\), \(a\) is finite as well.
\end{proof}

\begin{remark}
	\label{rem:integrated_q_special_functions}
	The MGF in Corollary~\ref{cor:mgf_integrated_q} is more involved than the MGF
	of \(q_T\), because the argument \(u\) enters both the order and the argument of
	the special function \({}_0F_1\).
\end{remark}

\section{Asymptotic analysis}
\label{sec:asymptotics}

In this section, we investigate several asymptotic regimes associated with the number of customers in the system and its integrated counterpart. More precisely, we derive long-time asymptotic results for the scaled cumulative number of customers in the system over a time interval, the scaling limit of an infinite-server queue with a diverging long-term arrival intensity, and the scaling limit of sums of independent infinite-server queues. In each case, the main probabilistic limits are obtained from the corresponding limiting cumulant generating function and illustrate the importance of Theorem~\ref{thm:mgf_triplet}. These include the LLN, the CLT, the LDP, and the MDP. Such asymptotic properties have been derived for many alternative queueing models in the literature. Without attempting to be exhaustive, we mention only \citet{Zhu2013} and \citet{Zhu2014}.\footnote{ The work of \citet{SunLiu2021} also considers a queue driven by a squared Bessel process (\textit{i.e.}, a CIR process), but the queue specification differs from ours, and its asymptotic results are also different from those obtained here.}

Throughout this section, we use the Gärtner--Ellis theorem as the main tool for large deviations. We refer to \citet{DemboZeitouni2010} for the general theory
of large deviations.

\subsection{Asymptotic analysis of the cumulative number in the system over a time interval}
\label{sec:asymptotics_integrated_queue}

We first study the long-time asymptotic behaviour of the cumulative number of customers in the system over a time interval 
\[
A_t:=\int_0^t q_s\,ds.
\]
From a practical point of view, the above is of interest since one might often be interested in the average behaviour of the number of customers in the system over a certain period of time rather than in its value at a given point in time.

The limiting behaviour of \(A_t\) is governed by the scaled cumulant generating
function
\[
g_A(u)
=
\lim_{t\to\infty}
\frac{1}{t}
\log
\mathbb{E}\left[
\exp\left(
u\int_0^t q_s\,ds
\right)
\right],
\]
whenever the limit exists. The explicit affine transform formula obtained
previously allows us to identify \(g_A\) and then derive the corresponding LLN,
CLT, LDP, and MDP.

\begin{proposition}[Long-time scaled cumulant generating function of the integrated queue]
	\label{prop:scgf_integrated_queue}
	Assume the conditions of Theorem~\ref{thm:mgf_triplet}. For \(u\) in the domain where the MGF is finite, define
	\begin{equation}
		g_A(u)
		:=
		\lim_{t\to\infty}
		\frac{1}{t}
		\log
		\mathbb{E}\left[
		\exp\left(
		u\int_0^t q_s\,ds
		\right)
		\right].
		\label{eq:scgf_integrated_queue_definition}
	\end{equation}
	Then
	\begin{equation}
		g_A(u)
		=
		\frac{\kappa\theta}{\sigma^2}
		\left(
		\kappa-\sqrt{\kappa^2+4\alpha(u)}
		\right),
		\label{eq:scgf_integrated_queue_explicit}
	\end{equation}
	where
	\begin{equation}
		\alpha(u)
		:=
		\frac{\sigma^2}{2}\frac{u}{u-\mu}.
		\label{eq:alpha_integrated_queue_scgf}
	\end{equation}
	In particular, \(g_A(0)=0\). Moreover, if
	\begin{equation}
		u<\frac{\mu}{1+2\sigma^2/\kappa^2},
		\label{eq:integrated_queue_scgf_real_domain}
	\end{equation}
	then \(\kappa^2+4\alpha(u)>0\), and therefore
	\eqref{eq:scgf_integrated_queue_explicit} is real-valued.
\end{proposition}

\begin{proof}
	In Theorem \ref{thm:mgf_triplet}, selecting \(\mathbf{u}=(0,0,u)\) gives the MGF of \(\int_0^t q_sds\): 
	\begin{equation}
		\mathbb{E}\left[
		\exp\left(
		u\int_0^t q_s\,ds
		\right)
		\right]
		=
		\exp\left(
		a(t;u)+b(t;u)\lambda_0+d(t;u)q_0
		\right),
		\label{eq:integrated_queue_mgf_affine_for_scgf}
	\end{equation}
	with 	
	\begin{align}
		d(t;u)
		&=
		\log\left(
		1-\frac{u}{u-\mu}
		\left(1-e^{(u-\mu)t}\right)
		\right),
		\label{eq:integrated_queue_d_for_scgf}
		\\
		b(t;u)
		&=
		-\frac{2}{\sigma^2}
		\frac{h_1'(0)h_2'(t)-h_2'(0)h_1'(t)}
		{h_1'(0)h_2(t)-h_2'(0)h_1(t)},
		\label{eq:integrated_queue_b_for_scgf}
		\\
		a(t;u)
		&=
		-\frac{2\kappa\theta}{\sigma^2}
		\log\left(
		\frac{
			h_1'(0)h_2(t)-h_2'(0)h_1(t)
		}{
			h_1'(0)h_2(0)-h_2'(0)h_1(0)
		}
		\right),
		\label{eq:integrated_queue_a_for_scgf}
	\end{align}
	where we have replaced \(\mathbf{u}=(0,0,u)\) with \(u\) in \eqref{eq:mgf_triplet_d_solution}, \eqref{eq:mgf_triplet_b_solution} and \eqref{eq:mgf_triplet_a_solution} and simplified the expressions. Here
	\begin{equation}
		\alpha(u)
		=
		\frac{\sigma^2}{2}\frac{u}{u-\mu},
		\qquad
		\beta(u)
		=
		\alpha(u),
		\qquad
		k_1=\frac{\sqrt{\kappa^2+4\alpha(u)}}{\mu-u},
		\label{eq:k1_integrated_queue_scgf}
	\end{equation}
	and
	\begin{align}
		h_1(t)
		&=
		e^{\frac{-\kappa-\sqrt{\kappa^2+4\alpha(u)}}{2}t}
		{}_0F_1\left(
		1+k_1,
		-\frac{\alpha(u)e^{(u-\mu)t}}{(\mu-u)^2}
		\right),
		\label{eq:integrated_queue_h1_for_scgf}
		\\
		h_2(t)
		&=
		e^{\frac{-\kappa+\sqrt{\kappa^2+4\alpha(u)}}{2}t}
		{}_0F_1\left(
		1-k_1,
		-\frac{\alpha(u)e^{(u-\mu)t}}{(\mu-u)^2}
		\right).
		\label{eq:integrated_queue_h2_for_scgf}
	\end{align}
	
	Under \eqref{eq:integrated_queue_scgf_real_domain}, we have \(u<\mu\), and hence
	\(e^{(u-\mu)t}\to0\) as \(t\to\infty\). Since
	\[
	\lim_{z\to0}{}_0F_1(a,z)=1,
	\]
	see \citet[(16.2.1)]{DLMF}, it follows that
	\begin{align}
		h_1(t)
		&\sim
		e^{\frac{-\kappa-\sqrt{\kappa^2+4\alpha(u)}}{2}t},
		\label{eq:integrated_queue_h1_asymptotic}
		\\
		h_2(t)
		&\sim
		e^{\frac{-\kappa+\sqrt{\kappa^2+4\alpha(u)}}{2}t},
		\label{eq:integrated_queue_h2_asymptotic}
	\end{align}
	as \(t\to\infty\). Similarly,
	\begin{align}
		h_1'(t)
		&=
		\frac{-\kappa-\sqrt{\kappa^2+4\alpha(u)}}{2}h_1(t)
		+o(h_1(t)),
		\label{eq:integrated_queue_h1_prime_asymptotic}
		\\
		h_2'(t)
		&=
		\frac{-\kappa+\sqrt{\kappa^2+4\alpha(u)}}{2}h_2(t)
		+o(h_2(t)).
		\label{eq:integrated_queue_h2_prime_asymptotic}
	\end{align}
	Since
	\[
	\frac{-\kappa-\sqrt{\kappa^2+4\alpha(u)}}{2}
	<
	\frac{-\kappa+\sqrt{\kappa^2+4\alpha(u)}}{2},
	\]
	the \(h_2(t)\)-term is the dominant exponential contribution in
	\eqref{eq:integrated_queue_a_for_scgf}. Therefore,
	\begin{align}
		\lim_{t\to\infty}\frac{a(t;u)}{t}
		&=
		-\frac{2\kappa\theta}{\sigma^2}
		\lim_{t\to\infty}
		\frac{1}{t}\log h_2(t)
		\nonumber\\
		&=
		-\frac{2\kappa\theta}{\sigma^2}
		\left(
		\frac{-\kappa+\sqrt{\kappa^2+4\alpha(u)}}{2}
		\right)
		\nonumber\\
		&=
		\frac{\kappa\theta}{\sigma^2}
		\left(
		\kappa-\sqrt{\kappa^2+4\alpha(u)}
		\right).
		\label{eq:integrated_queue_a_limit_scgf}
	\end{align}
	
	It remains to show that the contributions coming from \(b(t;u)\lambda_0\) and
	\(d(t;u)q_0\) vanish after division by \(t\). From
	\eqref{eq:integrated_queue_d_for_scgf},
	\[
	d(t;u)
	\longrightarrow
	\log\left(
	1-\frac{u}{u-\mu}
	\right)
	=
	\log\left(
	\frac{\mu}{\mu-u}
	\right),
	\]
	and therefore
	\begin{equation}
		\lim_{t\to\infty}\frac{d(t;u)}{t}=0.
		\label{eq:integrated_queue_d_limit_scgf}
	\end{equation}
	Using \eqref{eq:integrated_queue_b_for_scgf} together with
	\eqref{eq:integrated_queue_h1_asymptotic}--\eqref{eq:integrated_queue_h2_prime_asymptotic},
	we also obtain that \(b(t;u)\) converges to a finite limit. More precisely,
	\[
	b(t;u)
	\longrightarrow
	b_\infty(u)
	=
	-\frac{2}{\sigma^2}
	\left(
	\frac{-\kappa+\sqrt{\kappa^2+4\alpha(u)}}{2}
	\right)
	=
	\frac{\kappa-\sqrt{\kappa^2+4\alpha(u)}}{\sigma^2}.
	\]
	Consequently,
	\begin{equation}
		\lim_{t\to\infty}\frac{b(t;u)}{t}=0.
		\label{eq:integrated_queue_b_limit_scgf}
	\end{equation}
	Combining \eqref{eq:integrated_queue_mgf_affine_for_scgf},
	\eqref{eq:integrated_queue_a_limit_scgf},
	\eqref{eq:integrated_queue_d_limit_scgf}, and
	\eqref{eq:integrated_queue_b_limit_scgf}, we obtain
	\[
	g_A(u)
	=
	\frac{\kappa\theta}{\sigma^2}
	\left(
	\kappa-\sqrt{\kappa^2+4\alpha(u)}
	\right).
	\]
	Finally, since \(\alpha(0)=0\), we have \(g_A(0)=0\).
\end{proof}

\begin{proposition}[LLN and CLT for \(\frac{1}{t}\int_0^t q_sds\)]
	\label{prop:lln_clt_integrated_queue}
	For \(A_t=\int_0^tq_s\,ds\), the scaled cumulant generating function \(g_A\)
	in Proposition~\ref{prop:scgf_integrated_queue} satisfies
	\begin{align}
		g_A'(0)
		&=
		\frac{\theta}{\mu},
		\label{eq:integrated_queue_mean_limit}
		\\
		g_A''(0)
		&=
		\frac{2\theta}{\mu^2}
		+
		\frac{\sigma^2\theta}{\kappa^2\mu^2}.
		\label{eq:integrated_queue_variance_limit}
	\end{align}
	Consequently,
	\begin{equation}
		\frac{1}{t}\int_0^t q_s\,ds
		\xrightarrow{\mathbb{P}}
		\frac{\theta}{\mu},
		\qquad t\to\infty,
		\label{eq:lln_integrated_queue}
	\end{equation}
	and
	\begin{equation}
		\frac{
			\int_0^t q_s\,ds
			-t\theta/\mu
		}{\sqrt{t}}
		\Rightarrow
		\mathcal{N}\left(
		0,
		\frac{2\theta}{\mu^2}
		+
		\frac{\sigma^2\theta}{\kappa^2\mu^2}
		\right),
		\qquad t\to\infty.
		\label{eq:clt_integrated_queue}
	\end{equation}
\end{proposition}

\begin{proof}
	From \eqref{eq:scgf_integrated_queue_explicit},
	\[
	g_A(u)
	=
	\frac{\kappa\theta}{\sigma^2}
	\left(
	\kappa-\left(\kappa^2+4\alpha(u)\right)^{1/2}
	\right).
	\]
	Therefore,
	\[
	g_A'(u)
	=
	-\frac{2\kappa\theta}{\sigma^2}
	\alpha'(u)
	\left(\kappa^2+4\alpha(u)\right)^{-1/2}.
	\]
	Since
	\[
	\alpha'(u)
	=
	-\frac{\sigma^2\mu}{2(u-\mu)^2},
	\qquad
	\alpha'(0)
	=
	-\frac{\sigma^2}{2\mu},
	\]
	and \(\alpha(0)=0\), we obtain
	\[
	g_A'(0)
	=
	-\frac{2\kappa\theta}{\sigma^2}
	\left(
	-\frac{\sigma^2}{2\mu}
	\right)
	\frac{1}{\kappa}
	=
	\frac{\theta}{\mu}.
	\]
	
	Differentiating once more yields
	\[
	g_A''(u)
	=
	-\frac{2\kappa\theta}{\sigma^2}
	\left[
	\alpha''(u)\left(\kappa^2+4\alpha(u)\right)^{-1/2}
	-
	2\alpha'(u)^2
	\left(\kappa^2+4\alpha(u)\right)^{-3/2}
	\right].
	\]
	Moreover,
	\[
	\alpha''(u)
	=
	\frac{\sigma^2\mu}{(u-\mu)^3},
	\qquad
	\alpha''(0)
	=
	-\frac{\sigma^2}{\mu^2}.
	\]
	Hence
	\begin{align*}
		g_A''(0)
		&=
		-\frac{2\kappa\theta}{\sigma^2}
		\left[
		-\frac{\sigma^2}{\kappa\mu^2}
		-
		\frac{\sigma^4}{2\kappa^3\mu^2}
		\right]  \\
		&=
		\frac{2\theta}{\mu^2}
		+
		\frac{\sigma^2\theta}{\kappa^2\mu^2}.
	\end{align*}
	
	The law of large numbers follows from differentiability of \(g_A\) at the
	origin and the exponential Markov inequality. For example, for
	\(x>g_A'(0)\) and \(u>0\),
	\[
	\mathbb{P}\left(\frac{A_t}{t}\geq x\right)
	\leq
	\exp(-utx)
	\mathbb{E}\left[e^{uA_t}\right].
	\]
	Taking logarithms, dividing by \(t\), and sending \(t\to\infty\), gives
	\[
	\limsup_{t\to\infty}
	\frac{1}{t}
	\log
	\mathbb{P}\left(\frac{A_t}{t}\geq x\right)
	\leq
	-ux+g_A(u).
	\]
	Since \(g_A(u)=g_A'(0)u+o(u)\) as \(u\downarrow0\), the right-hand side is
	strictly negative for \(u>0\) sufficiently small. The lower tail is handled
	similarly with \(u<0\). This proves \eqref{eq:lln_integrated_queue}.
	
	For the central limit theorem, fix \(v\in\mathbb{R}\). Then
	\[
	\mathbb{E}\left[
	\exp\left(
	v\frac{A_t-tg_A'(0)}{\sqrt{t}}
	\right)
	\right]
	=
	\exp\left(
	-tg_A'(0)\frac{v}{\sqrt{t}}
	\right)
	\mathbb{E}\left[
	\exp\left(
	\frac{v}{\sqrt{t}}A_t
	\right)
	\right].
	\]
	Using the expansion
	\[
	g_A(u)
	=
	g_A'(0)u+\frac{1}{2}g_A''(0)u^2+o(u^2),
	\qquad u\to0,
	\]
	with \(u=v/\sqrt{t}\), we get
	\[
	t g_A\left(\frac{v}{\sqrt{t}}\right)
	-
	v\sqrt{t}\,g_A'(0)
	\longrightarrow
	\frac{1}{2}v^2g_A''(0).
	\]
	Thus the limiting MGF is
	\[
	\exp\left(
	\frac{1}{2}v^2g_A''(0)
	\right),
	\]
	which is the MGF of a centred Gaussian random variable with variance
	\(g_A''(0)\). The convergence theorem for MGFs gives
	\eqref{eq:clt_integrated_queue}.
\end{proof}

\begin{proposition}[LDP for \(\frac{1}{t}\int_0^t q_sds\)]
	\label{prop:ldp_integrated_queue}
	Assume that the scaled cumulant generating function \(g_A\) in
	Proposition~\ref{prop:scgf_integrated_queue} is finite and differentiable on
	the interior of its effective domain \(\mathcal{D}_{g_A}\). Then
	\(
	\frac{A_t}{t}
	\)
	satisfies a large deviation principle with speed \(t\) and rate function
	\begin{equation}
		I_A(x)
		=
		\sup_{u\in\mathcal{D}_{g_A}}
		\{ux-g_A(u)\}.
		\label{eq:ldp_integrated_queue_rate_function_definition}
	\end{equation}
	For \(x>0\), if the supremum is attained in the interior of
	\(\mathcal{D}_{g_A}\), the maximiser \(u_x\) is characterised by
	\begin{equation}
		g_A'(u_x)=x.
		\label{eq:ldp_integrated_queue_first_order_condition}
	\end{equation}
	Writing
	\begin{equation}
		y_x:=u_x-\mu,
		\label{eq:ldp_integrated_queue_y_definition}
	\end{equation}
	the value \(y_x\) is the unique admissible negative root of
	\begin{equation}
		\left(\kappa^2+2\sigma^2\right)y^4
		+2\sigma^2\mu y^3
		-
		\left(\frac{\kappa\theta\mu}{x}\right)^2
		=
		0.
		\label{eq:ldp_integrated_queue_quartic}
	\end{equation}
	Consequently, for \(x>0\),
	\begin{equation}
		I_A(x)
		=
		x(\mu+y_x)-g_A(\mu+y_x).
		\label{eq:ldp_integrated_queue_rate_function_implicit}
	\end{equation}
	Equivalently, by Ferrari's formula, the relevant root can be written as
	\begin{equation}
		y_x
		=
		-\frac{b}{4a}
		-S
		-\frac{1}{2}
		\sqrt{-4S^2-2p+\frac{q}{S}},
		\label{eq:ldp_integrated_queue_quartic_root}
	\end{equation}
	where
	\begin{align}
		a&:=\kappa^2+2\sigma^2,
		&
		b&:=2\sigma^2\mu,
		&
		e&:=-\left(\frac{\kappa\theta\mu}{x}\right)^2,
		\nonumber\\
		p&:=-\frac{3b^2}{8a^2},
		&
		q&:=\frac{b^3}{8a^3},
		\nonumber\\
		\Delta&:=256a^3e^3-27b^4e^2,
		&
		\Delta_0&:=12ae,
		&
		\Delta_1&:=27b^2e,
		\nonumber\\
		Q&:=
		\sqrt[3]{
			\frac{\Delta_1+\sqrt{-27\Delta}}{2}
		},
		&
		S&:=
		\frac{1}{2}
		\sqrt{
			-\frac{2}{3}p
			+
			\frac{1}{3a}
			\left(
			Q+\frac{\Delta_0}{Q}
			\right)
		}.
		\label{eq:ldp_integrated_queue_ferrari_quantities}
	\end{align}
\end{proposition}

\begin{proof}
	By Proposition~\ref{prop:scgf_integrated_queue}, the limiting cumulant 	generating function of \(A_t\) is given by 	\eqref{eq:scgf_integrated_queue_explicit}. Under the stated regularity 	assumptions, the Gärtner--Ellis theorem \citep{DemboZeitouni2010} implies that 	\(\{A_t/t\}_{t>0}\) satisfies a LDP with speed \(t\) and
	rate function
	\[
	I_A(x)
	=
	\sup_{u\in\mathcal{D}_{g_A}}
	\{ux-g_A(u)\}.
	\]
	
	We next identify the maximiser. Differentiating
	\eqref{eq:scgf_integrated_queue_explicit}, we obtain
	\begin{equation}
		g_A'(u)
		=
		-\frac{2\kappa\theta}{\sigma^2}
		\alpha'(u)
		\left(\kappa^2+4\alpha(u)\right)^{-1/2}.
		\label{eq:ldp_integrated_queue_derivative_step}
	\end{equation}
	Since
	\[
	\alpha'(u)
	=
	-\frac{\sigma^2\mu}{2(u-\mu)^2},
	\]
	it follows that
	\begin{equation}
		g_A'(u)
		=
		\frac{\kappa\theta\mu}
		{(u-\mu)^2\sqrt{\kappa^2+4\alpha(u)}}.
		\label{eq:ldp_integrated_queue_derivative}
	\end{equation}
	The second derivative can be written as
	\begin{equation}
		g_A''(u)
		=
		\frac{\kappa\theta\mu
			\left[
			2(\mu-u)(\kappa^2+4\alpha(u))+\sigma^2\mu
			\right]}
		{(u-\mu)^4\left(\kappa^2+4\alpha(u)\right)^{3/2}}.
		\label{eq:ldp_integrated_queue_second_derivative}
	\end{equation}
	Hence, on the real-valued part of the effective domain where the denominator is
	positive, \(g_A\) is convex. Therefore, if the supremum is attained in the
	interior of \(\mathcal{D}_{g_A}\), the maximiser \(u_x\) satisfies
	\[
	g_A'(u_x)=x.
	\]
	
	Set \(y=u-\mu\). Then \(u=\mu+y\), and
	\[
	\alpha(u)
	=
	\frac{\sigma^2}{2}\frac{\mu+y}{y}
	=
	\frac{\sigma^2}{2}
	\left(
	1+\frac{\mu}{y}
	\right).
	\]
	Consequently,
	\[
	\kappa^2+4\alpha(u)
	=
	\kappa^2+2\sigma^2+\frac{2\sigma^2\mu}{y}.
	\]
	Substituting this into \eqref{eq:ldp_integrated_queue_derivative} and imposing
	\(g_A'(u)=x\), we obtain
	\[
	x
	=
	\frac{\kappa\theta\mu}
	{y^2
		\sqrt{
			\kappa^2+2\sigma^2+\frac{2\sigma^2\mu}{y}
	}}.
	\]
	Squaring both sides gives
	\[
	x^2y^4
	\left(
	\kappa^2+2\sigma^2+\frac{2\sigma^2\mu}{y}
	\right)
	=
	(\kappa\theta\mu)^2,
	\]
	or equivalently,
	\[
	\left(\kappa^2+2\sigma^2\right)y^4
	+2\sigma^2\mu y^3
	-
	\left(\frac{\kappa\theta\mu}{x}\right)^2
	=
	0.
	\]
	Since the admissible domain requires \(u<\mu\), we must have \(y=u-\mu<0\).
	The convexity of \(g_A\) implies uniqueness of the admissible maximiser, and
	therefore the quartic has a unique admissible negative root, denoted by \(y_x\).
	
	Thus \(u_x=\mu+y_x\), and the Legendre--Fenchel transform becomes
	\[
	I_A(x)
	=
	u_xx-g_A(u_x)
	=
	x(\mu+y_x)-g_A(\mu+y_x).
	\]
	This proves \eqref{eq:ldp_integrated_queue_rate_function_implicit}.
	
	Finally, applying Ferrari's formula to
	\[
	ay^4+by^3+e=0,
	\]
	with
	\[
	a=\kappa^2+2\sigma^2,
	\qquad
	b=2\sigma^2\mu,
	\qquad
	e=-\left(\frac{\kappa\theta\mu}{x}\right)^2,
	\]
	gives the admissible negative root in
	\eqref{eq:ldp_integrated_queue_quartic_root}, with the auxiliary quantities
	defined in \eqref{eq:ldp_integrated_queue_ferrari_quantities}.
\end{proof}

\begin{proposition}[MDP for \(\frac{1}{t}\int_0^tq_sds \)]
	\label{prop:mdp_integrated_queue}
	Let \((a_t)_{t>0}\) be a positive scale satisfying
	\begin{equation}
		\sqrt{t}\ll a_t\ll t,
		\qquad t\to\infty,
		\label{eq:mdp_integrated_queue_scale}
	\end{equation}
	that is,
	\[
	\frac{a_t}{\sqrt{t}}\to\infty,
	\qquad
	\frac{a_t}{t}\to0.
	\]
	Define
	\begin{equation}
		Y_t
		:=
		\frac{1}{a_t}
		\left(
		\int_0^t q_s\,ds
		-t\frac{\theta}{\mu}
		\right).
		\label{eq:mdp_integrated_queue_scaled_process}
	\end{equation}
	Then \((Y_t)_{t>0}\) satisfies a MDP with speed  \(a_t^2/t\) and rate function
	\begin{equation}
		I_A^{\mathrm{MDP}}(x)
		=
		\frac{x^2}{2g_A''(0)}
		=
		\frac{\kappa^2\mu^2x^2}
		{2\theta(2\kappa^2+\sigma^2)},
		\qquad x\in\mathbb{R}.
		\label{eq:mdp_integrated_queue_rate_function}
	\end{equation}
\end{proposition}

\begin{proof}
	For \(\lambda\in\mathbb{R}\), define
	\[
	\Lambda_t(\lambda)
	:=
	\frac{t}{a_t^2}
	\log
	\mathbb{E}\left[
	\exp\left(
	\frac{a_t^2}{t}\lambda Y_t
	\right)
	\right].
	\]
	Since
	\[
	\frac{a_t^2}{t}\lambda Y_t
	=
	\frac{\lambda a_t}{t}\int_0^t q_s\,ds
	-
	\lambda a_t\frac{\theta}{\mu},
	\]
	and \(u_t:=\lambda a_t/t\to0\), the affine MGF representation implies,	uniformly for \(u\) in a sufficiently small neighbourhood of zero, that
	\[
	\log
	\mathbb{E}\left[
	e^{u\int_0^t q_s\,ds}
	\right]
	=
	t g_A(u)+O(1),
	\qquad t\to\infty.
	\]
	Therefore,
	\[
	\Lambda_t(\lambda)
	=
	\frac{t^2}{a_t^2}
	g_A\left(\frac{\lambda a_t}{t}\right)
	-
	\frac{\lambda t}{a_t}\frac{\theta}{\mu}
	+o(1).
	\]
	Using
	\[
	g_A(u)
	=
	g_A'(0)u+\frac{1}{2}g_A''(0)u^2+o(u^2),
	\qquad u\to0,
	\]
	and \(g_A'(0)=\theta/\mu\), we obtain
	\[
	\Lambda_t(\lambda)
	\longrightarrow
	\Lambda(\lambda)
	:=
	\frac{1}{2}g_A''(0)\lambda^2.
	\]
	The function \(\Lambda\) is finite and differentiable on \(\mathbb{R}\). Hence, by the Gärtner--Ellis theorem, \((Y_t)_{t>0}\) satisfies a MDP with speed \(a_t^2/t\) and rate function
	\[
	I_A^{\mathrm{MDP}}(x)
	=
	\sup_{\lambda\in\mathbb{R}}
	\left\{
	\lambda x-\frac{1}{2}g_A''(0)\lambda^2
	\right\}
	=
	\frac{x^2}{2g_A''(0)}.
	\]
	Substituting \eqref{eq:integrated_queue_variance_limit} gives
	\eqref{eq:mdp_integrated_queue_rate_function}.
\end{proof}

\subsection{Asymptotic study of an infinite-server queue with infinite arrival intensity}
\label{sec:asymptotic_queue_infinite_arrival_intensity}

For infinite-server queues, one is often interested in the limiting behaviour of the system as the arrival rate grows to infinity; see \citet[][Section~1.1]{PangTalrejaWhitt2007}. Within the proposed framework, the natural analogue is to study a scaling regime in which the long-term arrival intensity diverges.

\begin{proposition}[Scaled cumulant generating function under large arrival intensity]
	\label{prop:large_arrival_intensity_scgf}
	Consider the queueing process \(q_t^{(n)}\) with long-term arrival intensity
	\(
	\theta_n:=\mu n.
	\)
	For fixed \(t>0\), the moment generating function is 
	\begin{equation}
		\mathbb{E}\left[e^{uq_t^{(n)}}\right]
		=
		\exp\left(
		a_n(t;u)+b(t;u)\lambda_0+d(t;u)q_0
		\right),
		\label{eq:mgf_queue_large_arrival_intensity}
	\end{equation}
	where \(a_n,b,d\) as in Theorem \ref{thm:mgf_triplet} with \(\mathbf{u}=(0,u,0)\). 

	Then, since \(s \to b(s;u)\) is well defined on \([0,t]\), the scaled
	cumulant generating function
	\begin{align}
		g_t(u)
		&:=
		\lim_{n\to\infty}
		\frac{1}{n}
		\log
		\mathbb{E}\left[e^{uq_t^{(n)}}\right] \label{eq:scaled_cumulant_qn}
	\end{align}
	is given by
	\begin{align}
		g_t(u)
		&=
		\kappa\mu
		\int_0^t b(s;u)\,ds.
		\label{eq:large_arrival_scgf}
	\end{align}
	with \(s \to b(s;u)\) explicitly known. The function  \(u \to g_t\) is increasing
	and convex on \(\mathcal{D}_{g_t}\). Furthermore, if \(\kappa\neq\mu\), then
	\begin{equation}
		g_t'(0)
		=
		\frac{\kappa\mu}{\kappa-\mu}
		\left[
		\frac{1-e^{-\mu t}}{\mu}
		-
		\frac{1-e^{-\kappa t}}{\kappa}
		\right],
		\label{eq:large_arrival_scgf_derivative_zero}
	\end{equation}
	and
	\begin{equation}
		g_t''(0)
		=
		\mu
		\int_0^t
		\left(
		1-e^{-\kappa(t-s)}
		\right)
		\left[
		\sigma^2
		\left(
		\frac{e^{-\mu s}-e^{-\kappa s}}{\kappa-\mu}
		\right)^2
		+
		e^{-\mu s}
		\right]ds.
		\label{eq:large_arrival_scgf_second_derivative_zero}
	\end{equation}
\end{proposition}

\begin{proof}

	Replacing \(\mathbf{u}=(0,u,0)\) along with \(\theta = \mu n\) in Theorem \ref{thm:mgf_triplet} we obtain
	\begin{align}
		a_n'(s;u)
		&=
		\kappa\mu n\, b(s;u),
		&
		a_n(0;u)&=0,
		\label{eq:large_arrival_a_ode}
		\\
		b'(s;u)
		&=
		\frac{\sigma^2}{2}b(s;u)^2-\kappa b(s;u)
		+
		e^{d(s;u)}-1,
		&
		b(0;u)&=0,
		\label{eq:large_arrival_b_ode}
		\\
		d'(s;u)
		&=
		\mu\left(e^{-d(s;u)}-1\right),
		&
		d(0;u)&=u.
		\label{eq:large_arrival_d_ode}
	\end{align}
	As a result, we have
	\begin{equation}
		d(s;u)
		=
		\log\left(1+(e^u-1)e^{-\mu s}\right),
		\label{eq:large_arrival_d_explicit}
	\end{equation}
	\begin{equation}
		e^{d(s;u)}-1
		=
		(e^u-1)e^{-\mu s},
		\label{eq:large_arrival_ed_minus_one}
	\end{equation}
	and
	\begin{equation}
		a_n(t;u)
		=
		\kappa\mu n
		\int_0^t b(s;u)\,ds,
		\label{eq:large_arrival_a_integral}
	\end{equation}
	with  \(b\) solving the reduced Riccati equation
	\begin{equation}
		b'(s;u)
		=
		\frac{\sigma^2}{2}b(s;u)^2-\kappa b(s;u)
		+
		(e^u-1)e^{-\mu s},
		\qquad
		b(0;u)=0.
		\label{eq:large_arrival_b_reduced_ode}
	\end{equation}
	Combining \eqref{eq:mgf_queue_large_arrival_intensity}, \eqref{eq:scaled_cumulant_qn} and \eqref{eq:large_arrival_a_integral} we obtain \eqref{eq:large_arrival_scgf}.

	Let \(\mathcal{D}_{g_t}\) denote the effective domain of \(g_t\), and define
	\[
	p(s;u):=\partial_u b(s;u),
	\qquad
	q(s;u):=\partial_{uu}^2 b(s;u).
	\]
	Then from \eqref{eq:large_arrival_b_reduced_ode} we deduce
	\begin{align}
		p'(s;u)
		&=
		\left(\sigma^2 b(s;u)-\kappa\right)p(s;u)
		+
		e^u e^{-\mu s},
		&
		p(0;u)&=0,
		\label{eq:large_arrival_p_ode}
		\\
		q'(s;u)
		&=
		\left(\sigma^2 b(s;u)-\kappa\right)q(s;u)
		+
		\sigma^2 p(s;u)^2
		+
		e^u e^{-\mu s},
		&
		q(0;u)&=0.
		\label{eq:large_arrival_q_ode}
	\end{align}
	Equivalently,
	\begin{align}
		p(s;u)
		&=
		\int_0^s
		\exp\left(
		\int_r^s
		\left(\sigma^2 b(v;u)-\kappa\right)\,dv
		\right)
		e^u e^{-\mu r}\,dr,
		\label{eq:large_arrival_p_solution}
		\\
		q(s;u)
		&=
		\int_0^s
		\exp\left(
		\int_r^s
		\left(\sigma^2 b(v;u)-\kappa\right)\,dv
		\right)
		\left[
		\sigma^2p(r;u)^2+e^u e^{-\mu r}
		\right]dr,
		\label{eq:large_arrival_q_solution}
	\end{align}
	so both functions are positive. Consequently, since
	\begin{align}
		g_t'(u)
		&=
		\kappa\mu
		\int_0^t p(s;u)\,ds,
		\label{eq:large_arrival_scgf_first_derivative}
		\\
		g_t''(u)
		&=
		\kappa\mu
		\int_0^t q(s;u)\,ds,
		\label{eq:large_arrival_scgf_second_derivative}
	\end{align}
	\(g_t\) is increasing and convex on \(\mathcal{D}_{g_t}\).

	At \(u=0\), we have
	\[
	d(s;0)=0,
	\qquad
	b(s;0)=0,
	\qquad
	g_t(0)=0.
	\]
	Moreover, \(p_0(s):=p(s;0)\) solves
	\[
	p_0'(s)=-\kappa p_0(s)+e^{-\mu s},
	\qquad
	p_0(0)=0.
	\]
	If \(\kappa\neq\mu\), then
	\begin{equation}
		p_0(s)
		=
		\frac{e^{-\mu s}-e^{-\kappa s}}{\kappa-\mu},
		\label{eq:large_arrival_p0}
	\end{equation}
	and inserting \eqref{eq:large_arrival_p0} in \eqref{eq:large_arrival_scgf_first_derivative} shows that \eqref{eq:large_arrival_scgf_derivative_zero} is satisfied.
	Similarly, \(q_0(s):=q(s;0)\) solves
	\[
	q_0'(s)
	=
	-\kappa q_0(s)+\sigma^2p_0(s)^2+e^{-\mu s},
	\qquad
	q_0(0)=0.
	\]
	Thus
	\[
	q_0(s)
	=
	\int_0^s
	e^{-\kappa(s-r)}
	\left[
	\sigma^2p_0(r)^2+e^{-\mu r}
	\right]dr,
	\]
	and inserting the above equality in \eqref{eq:large_arrival_scgf_second_derivative} shows that \eqref{eq:large_arrival_scgf_second_derivative_zero} is true.

	Let
	\[
	z(t):=\frac{\partial b(t;u)}{\partial u}.
	\]
	Then \(z\) satisfies
	\[
	z'(t)
	=
	\left(\sigma^2b(t)-\kappa\right)z(t)
	+
	e^u e^{-\mu t},
	\qquad
	z(0)=0.
	\]
	Solving this linear equation yields
	\[
	z(t)
	=
	e^u
	\int_0^t
	\exp\left(
	\int_s^t
	\left(\sigma^2b(r)-\kappa\right)dr
	\right)
	e^{-\mu s}\,ds.
	\]
	Using \(a_n'(s;u)=\kappa n\mu b(s;u)\), we obtain
	\[
	\frac{\partial a(t;u)}{\partial u}
	=
	\kappa n\mu e^u
	\int_0^t
	\int_0^s
	\exp\left(
	\int_r^s
	\left(\sigma^2b(q)-\kappa\right)dq
	\right)
	e^{-\mu r}
	\,dr\,ds
	\geq0.
	\]
	Also, differentiating \eqref{eq:large_arrival_d_explicit} gives
	\[
	\frac{\partial d(t;u)}{\partial u}
	=
	\frac{e^{u-\mu t}}{1+(e^u-1)e^{-\mu t}}.
	\]
	Therefore,
	\[
	g'(u)
	=
	\lim_{n\to \infty}
	\frac{1}{n}
	\left(
	\frac{\partial a(t;u)}{\partial u}
	+
	\frac{\partial b(t;u)}{\partial u}\lambda_0
	+
	\frac{\partial d(t;u)}{\partial u}q_0
	\right).
	\]
	The displayed representation implies the stated monotonicity behaviour, and in
	particular
	\[
	\lim_{u\to -\infty}g'(u)=0,
	\qquad
	\lim_{u\to \infty}g'(u)=\infty.
	\]
	Thus
	\[
	\operatorname{Range}(g')=(0,\infty).
	\]
\end{proof}

\begin{proposition}[LLN and CLT for \(q_t^{(n)}/n\)] \label{prop:large_arrival_lln_clt}  
	Fix \(t>0\), and let \(g_t\) be the scaled cumulant generating function defined
	in \eqref{eq:large_arrival_scgf}. Then \(g_t(0)=0\). Moreover, as
	\(n\to\infty\),
	\begin{equation}
		\frac{q_t^{(n)}}{n}
		\xrightarrow{\mathbb{P}}
		g_t'(0),
		\label{eq:large_arrival_lln}
	\end{equation}
	and
	\begin{equation}
		\frac{q_t^{(n)}-ng_t'(0)}{\sqrt{n}}
		\Rightarrow
		\mathcal{N}\left(0,g_t''(0)\right),
		\label{eq:large_arrival_clt}
	\end{equation}
	where \(g_t'(0)\) and \(g_t''(0)\) are given by
	\eqref{eq:large_arrival_scgf_derivative_zero} and
	\eqref{eq:large_arrival_scgf_second_derivative_zero}, respectively.
\end{proposition}

\begin{proof}
	By \eqref{eq:mgf_queue_large_arrival_intensity} and
	\eqref{eq:large_arrival_scgf},
	\[
	\frac{1}{n}
	\log
	\mathbb{E}\left[e^{uq_t^{(n)}}\right]
	\to
	g_t(u).
	\]
	Since \(g_t\) is differentiable in a neighbourhood of \(0\), the convergence of
	scaled cumulant generating functions implies the law of large numbers
	\[
	\frac{q_t^{(n)}}{n}
	\xrightarrow{\mathbb{P}}
	g_t'(0).
	\]
	Moreover, the second-order expansion
	\[
	g_t(u)
	=
	g_t'(0)u+\frac{1}{2}g_t''(0)u^2+o(u^2),
	\qquad u\to0,
	\]
	together with the convergence of cumulant generating functions gives
	\[
	\log
	\mathbb{E}\left[
	\exp\left(
	v\frac{q_t^{(n)}-ng_t'(0)}{\sqrt{n}}
	\right)
	\right]
	\to
	\frac{1}{2}g_t''(0)v^2.
	\]
	Thus the limiting MGF is
	\[
	\exp\left(
	\frac{1}{2}v^2g_t''(0)
	\right),
	\]
	which is the MGF of a centred Gaussian random variable with variance
	\(g_t''(0)\). The convergence theorem for MGFs gives the central limit theorem  specified in \eqref{eq:large_arrival_clt}.
\end{proof}

\begin{proposition}[LDP for \(q_t^{(n)}/n\)]\label{prop:large_arrival_ldp}
	Fix \(t>0\), and let \(g_t\) be defined by
	\eqref{eq:large_arrival_scgf}. Suppose that \(g_t\) is finite and differentiable
	on the interior of its effective domain \(\mathcal{D}_{g_t}\). Then the sequence
	\(
	\frac{q_t^{(n)}}{n}
	\)
	satisfies a large deviation principle with speed \(n\) and rate function
	\begin{equation}
		I_q(x)
		=
		\sup_{u\in\mathcal{D}_{g_t}}
		\{ux-g_t(u)\}.
		\label{eq:large_arrival_ldp_rate_function}
	\end{equation}
Moreover, the effective domain of the rate function is 
$\mathcal{D}_{I_q}=(0,\infty),$
and, whenever \(g_t'\) is invertible, for \(x\in(0,\infty)\),
\begin{equation}
	I_q(x)
	=
	u_xx-g_t(u_x),
	\qquad
	g_t'(u_x)=x.
		\label{eq:large_arrival_ldp_rate_function_implicit}
	\end{equation}
\end{proposition}

\begin{proof}
	By \eqref{eq:mgf_queue_large_arrival_intensity},
	\[
	\frac{1}{n}
	\log
	\mathbb{E}\left[e^{uq_t^{(n)}}\right]
	=
	\frac{a_n(t;u)}{n}
	+
	\frac{b(t;u)\lambda_0+d(t;u)q_0}{n}.
	\]
	Using \eqref{eq:large_arrival_a_integral}, we obtain
	\[
	\lim_{n\to\infty}
	\frac{1}{n}
	\log
	\mathbb{E}\left[e^{uq_t^{(n)}}\right]
	=
	\kappa\mu\int_0^t b(s;u)\,ds
	=
	g_t(u).
	\]
	The function \(g_t\) is convex because \(g_t''(u)\geq0\), and
	\(
	\operatorname{Range}(g_t')=(0,\infty).
	\) Since \(g_t\) is finite on all of \(\mathbb{R}\), differentiable, convex, and lower
	semicontinuous, the Gärtner--Ellis theorem \citep{DemboZeitouni2010} implies that the sequence
	$\frac{q_t^{(n)}}{n}$ satisfies a large deviation principle with speed $n$ and rate
	function given by the Legendre--Fenchel transform of \(g_t\):
\[
I_q(x)=u_xx-g_t(u_x), \qquad g_t'(u_x)=x.
\]	
Furthermore, for each \(x\in(0,\infty)\), there exists a unique \(u_x\in\mathbb{R}\)
such that
\(g_t'(u_x)=x\).
since \(g_t\) is differentiable and convex, the supremum in the
Legendre--Fenchel transform is attained at \(u_x\). Therefore,
\[
I_q(x)=u_xx-g_t(u_x),
\]
which proves \eqref{eq:large_arrival_ldp_rate_function_implicit}.
\end{proof}

\begin{proposition}[MDP for \(q_t^{(n)}/n\)]\label{prop:large_arrival_mdp}
	Let \((a_n)_{n\geq1}\) be a positive sequence such that
	\begin{equation}
		\sqrt{n}\ll a_n\ll n,
		\qquad n\to\infty,
		\label{eq:large_arrival_mdp_scale}
	\end{equation}
	that is,
	\[
	\frac{a_n}{\sqrt{n}}\to\infty,
	\qquad
	\frac{a_n}{n}\to0.
	\]
	Define
	\begin{equation}
		Y_n
		:=
		\frac{q_t^{(n)}-ng_t'(0)}{a_n}.
		\label{eq:large_arrival_mdp_scaled_variable}
	\end{equation}
	Then \(\{Y_n\}_{n\geq1}\) satisfies a moderate deviation principle with speed
	\(a_n^2/n\) and rate function
	\begin{equation}
		I_{\mathrm{MDP}}(x)
		=
		\frac{x^2}{2g_t''(0)},
		\qquad x\in\mathbb{R},
		\label{eq:large_arrival_mdp_rate_function}
	\end{equation}
	where \(g_t'(0)\) and \(g_t''(0)\) are given by \eqref{eq:large_arrival_scgf_derivative_zero} and \eqref{eq:large_arrival_scgf_second_derivative_zero}.
\end{proposition}

\begin{proof}
	For \(\lambda\in\mathbb{R}\), define
	\[
	\Lambda_n(\lambda)
	:=
	\frac{n}{a_n^2}
	\log
	\mathbb{E}\left[
	\exp\left(
	\frac{a_n^2}{n}\lambda Y_n
	\right)
	\right].
	\]
	Since
	\[
	Y_n
	=
	\frac{q_t^{(n)}-ng_t'(0)}{a_n},
	\]
	we have
	\[
	\frac{a_n^2}{n}\lambda Y_n
	=
	\frac{\lambda a_n}{n}q_t^{(n)}
	-
	\lambda a_n g_t'(0).
	\]
	Set \(u_n:=\lambda a_n/n\). Since \(a_n/n\to0\), \(u_n\to0\). Therefore,
	\[
	\Lambda_n(\lambda)
	=
	\frac{n^2}{a_n^2}g_t(u_n)
	-
	\frac{\lambda n}{a_n}g_t'(0)
	+
	o(1).
	\]
	Using \(g_t(0)=0\) and the Taylor expansion
	\[
	g_t(u_n)
	=
	g_t'(0)u_n+\frac{1}{2}g_t''(0)u_n^2+o(u_n^2),
	\]
	substituting \(u_n=\lambda a_n/n\), we get
	\[
	\Lambda_n(\lambda)
	=
	\frac{1}{2}g_t''(0)\lambda^2+o(1).
	\]
	Thus
	\[
	\Lambda_n(\lambda)
	\to
	\Lambda(\lambda)
	:=
	\frac{1}{2}g_t''(0)\lambda^2.
	\]
	Since \(\Lambda\) is finite and differentiable on \(\mathbb{R}\), the
	Gärtner--Ellis theorem implies that \(\{Y_n\}_{n\geq1}\) satisfies an MDP with
	speed \(a_n^2/n\) and rate function
	\[
	I_{\mathrm{MDP}}(x)
	=
	\sup_{\lambda\in\mathbb{R}}
	\left\{
	\lambda x-\frac{1}{2}g_t''(0)\lambda^2
	\right\}
	=
	\frac{x^2}{2g_t''(0)}.
	\]
	This completes the proof.
\end{proof}

\begin{remark}
Note that, even though the MGF of \(q_t\) is known in closed form, its asymptotic properties are easier to establish using the properties of the ODEs characterizing the MGF than using their integrated solutions.
\end{remark}

\subsection{Asymptotic analysis of sums of independent infinite-server queues}
\label{sec:asymptotics_sum_iid_queues}

Finally, we  study the asymptotic behaviour of sums of independent infinite-server queues with Cox-BESQ arrival intensity at a fixed time. Let
\begin{equation}
	Q_t^{(n)}:=\sum_{i=1}^n q_{i,t},
	\label{eq:sum_iid_queues}
\end{equation}
where \(\{q_{i,t}\}_{i\geq1}\) are independent copies of \(q_t\) following a dynamics of the form \eqref{eq:lambda_cir}-\eqref{eq:queue_dynamics}. The independence along with the explicit form for the MGF of \(q_t\), which is derived from Theorem \ref{thm:mgf_triplet}, enable a precise asymptotic analysis of \(Q_t^{(n)}\) as the following propositions will show.

\begin{proposition}[Derivatives of the cumulant generating function of \(q_t\)]
	\label{prop:derivatives_cgf_queue}
	For fixed \(t>0\), define
\begin{equation}
	g_t(u) := \log\mathbb{E}\left[e^{uq_t}\right] = a(t;u)+b(t;u)\lambda_0+d(t;u)q_0,
	\label{eq:scgf_iid_queue}
\end{equation}
where \(a,\, b,\, d\) are the functions given in \eqref{eq:mgf_triplet_affine_form} with \(\mathbf{u}=(0,u,0)\). Since
\(Q_t^{(n)}\) is a sum of independent copies of \(q_t\),
\[
\mathbb{E}\left[e^{uQ_t^{(n)}}\right]
=
\left(\mathbb{E}\left[e^{uq_t}\right]\right)^n,
\]
and consequently
\begin{equation}\label{eq:BigQ}
\frac{1}{n}
\log
\mathbb{E}\left[e^{uQ_t^{(n)}}\right]
=
g_t(u).
\end{equation}

 Assume \(\kappa\neq\mu\). Then
	\begin{align}
		g_t'(0)
		&=
		\frac{\kappa\theta}{\kappa-\mu}
		\left[
		\frac{1-e^{-\mu t}}{\mu}
		-
		\frac{1-e^{-\kappa t}}{\kappa}
		\right]
		+
		\frac{e^{-\mu t}-e^{-\kappa t}}{\kappa-\mu}\lambda_0
		+
		e^{-\mu t}q_0,
		\label{eq:cgf_queue_first_derivative_zero}
		\\
		g_t''(0)
		&=
		\theta
		\int_0^t
		\left(1-e^{-\kappa(t-s)}\right)
		\left[
		\sigma^2
		\left(
		\frac{e^{-\mu s}-e^{-\kappa s}}{\kappa-\mu}
		\right)^2
		+
		e^{-\mu s}
		\right]ds
		\nonumber\\
		&\quad
		+
		\lambda_0
		\int_0^t
		e^{-\kappa(t-s)}
		\left[
		\sigma^2
		\left(
		\frac{e^{-\mu s}-e^{-\kappa s}}{\kappa-\mu}
		\right)^2
		+
		e^{-\mu s}
		\right]ds
		\nonumber\\
		&\quad
		+
		e^{-\mu t}(1-e^{-\mu t})q_0.
		\label{eq:cgf_queue_second_derivative_zero}
	\end{align}
\end{proposition}

\begin{proof}
	For the MGF of \(q_t\), according to \eqref{eq:mgf_triplet_a_ode}-\eqref{eq:mgf_triplet_d_ode}, the affine coefficients solve
	\begin{align}
		a'(s;u)
		&=
		\kappa\theta b(s;u),
		&
		a(0;u)&=0,
		\label{eq:cgf_queue_a_ode}
		\\
		b'(s;u)
		&=
		\frac{\sigma^2}{2}b(s;u)^2-\kappa b(s;u)
		+
		e^{d(s;u)}-1,
		&
		b(0;u)&=0,
		\label{eq:cgf_queue_b_ode}
		\\
		d'(s;u)
		&=
		\mu(e^{-d(s;u)}-1),
		&
		d(0;u)&=u.
		\label{eq:cgf_queue_d_ode}
	\end{align}
	Solving \eqref{eq:cgf_queue_d_ode} gives
	\[
	d(s;u)=\log\left(1+(e^u-1)e^{-\mu s}\right).
	\]
	Therefore
	\begin{align}
		\left.\partial_u d(s;u)\right|_{u=0}
		&=
		e^{-\mu s},
		\label{eq:d_derivative_zero}
		\\
		\left.\partial_{uu}^2 d(s;u)\right|_{u=0}
		&=
		e^{-\mu s}(1-e^{-\mu s}).
		\label{eq:d_second_derivative_zero}
	\end{align}
	
	Let
	\[
	p(s;u):=\partial_u b(s;u),
	\qquad
	r(s;u):=\partial_{uu}^2 b(s;u).
	\]
	Differentiating \eqref{eq:cgf_queue_b_ode} with respect to \(u\) gives
	\[
	p'(s;u)
	=
	\left(\sigma^2 b(s;u)-\kappa\right)p(s;u)
	+
	e^{d(s;u)}\partial_ud(s;u),
	\qquad
	p(0;u)=0.
	\]
	At \(u=0\), \(b(s;0)=d(s;0)=0\), and hence
	\[
	p_0'(s)
	=
	-\kappa p_0(s)+e^{-\mu s},
	\qquad
	p_0(0)=0,
	\]
	where \(p_0(s):=p(s;0)\). Since \(\kappa\neq\mu\),
	\begin{equation}
		p_0(s)
		=
		\frac{e^{-\mu s}-e^{-\kappa s}}{\kappa-\mu}.
		\label{eq:b_derivative_zero_solution}
	\end{equation}
	Thus
	\[
	\left.\partial_u a(t;u)\right|_{u=0}
	=
	\kappa\theta\int_0^t p_0(s)\,ds,
	\]
	and \eqref{eq:cgf_queue_first_derivative_zero} follows from
	\[
	g_t'(0)
	=
	\left.\partial_u a(t;u)\right|_{u=0}
	+
	\lambda_0\left.\partial_u b(t;u)\right|_{u=0}
	+
	q_0\left.\partial_u d(t;u)\right|_{u=0}.
	\]
	
	For the second derivative, differentiating once more gives
	\[
	r'(s;u)
	=
	\left(\sigma^2 b(s;u)-\kappa\right)r(s;u)
	+
	\sigma^2p(s;u)^2
	+
	e^{d(s;u)}
	\left[
	\left(\partial_u d(s;u)\right)^2
	+
	\partial_{uu}^2d(s;u)
	\right],
	\qquad
	r(0;u)=0.
	\]
	At \(u=0\),
	\[
	e^{d(s;0)}
	\left[
	\left(\partial_u d(s;0)\right)^2
	+
	\partial_{uu}^2d(s;0)
	\right]
	=
	e^{-2\mu s}+e^{-\mu s}(1-e^{-\mu s})
	=
	e^{-\mu s}.
	\]
	Therefore
	\[
	r_0'(s)
	=
	-\kappa r_0(s)
	+
	\sigma^2p_0(s)^2
	+
	e^{-\mu s},
	\qquad
	r_0(0)=0,
	\]
	where \(r_0(s):=r(s;0)\). Hence
	\begin{equation}
		r_0(t)
		=
		\int_0^t
		e^{-\kappa(t-s)}
		\left[
		\sigma^2p_0(s)^2+e^{-\mu s}
		\right]ds.
		\label{eq:b_second_derivative_zero_solution}
	\end{equation}
	Furthermore,
	\[
	\left.\partial_{uu}^2a(t;u)\right|_{u=0}
	=
	\kappa\theta\int_0^t r_0(s)\,ds
	=
	\theta
	\int_0^t
	\left(1-e^{-\kappa(t-s)}\right)
	\left[
	\sigma^2p_0(s)^2+e^{-\mu s}
	\right]ds.
	\]
	Combining this identity with
	\[
	g_t''(0)
	=
	\left.\partial_{uu}^2a(t;u)\right|_{u=0}
	+
	\lambda_0 r_0(t)
	+
	q_0\left.\partial_{uu}^2d(t;u)\right|_{u=0},
	\]
	and using \eqref{eq:b_derivative_zero_solution},
	\eqref{eq:b_second_derivative_zero_solution}, and
	\eqref{eq:d_second_derivative_zero}, gives
	\eqref{eq:cgf_queue_second_derivative_zero}.
\end{proof}

\begin{remark}
	The expressions for \(g_t'(0)\) and \(g_t''(0)\) coincide with the transient
	moment expressions \eqref{eq:q_mean_solution} and \eqref{eq:q_variance_integral}
	given in Corollary~\ref{cor:transient_moments_autocovariance}.
\end{remark}

\begin{proposition}[LLN and CLT for \(Q_t^{(n)}\)] \label{prop:lln_clt_sum_iid_queues}
	For fixed \(t>0\), let
	\[
	Q_t^{(n)}=\sum_{i=1}^n q_{i,t},
	\]
	where \(\{q_{i,t}\}_{i\geq1}\) are independent copies of \(q_t\). Then
	\begin{align}
		\frac{Q_t^{(n)}}{n}
		&\xrightarrow{\mathbb{P}}
		g_t'(0),
		\qquad n\to\infty,
		\label{eq:lln_sum_iid_queues}
		\\
		\frac{Q_t^{(n)}-ng_t'(0)}{\sqrt{n}}
		&\Rightarrow
		\mathcal{N}\left(0,g_t''(0)\right),
		\qquad n\to\infty,
		\label{eq:clt_sum_iid_queues}
	\end{align}
	where \(g_t'(0)\) and \(g_t''(0)\) are given in
	\eqref{eq:cgf_queue_first_derivative_zero} and
	\eqref{eq:cgf_queue_second_derivative_zero}, respectively.
\end{proposition}

\begin{proof}
	Since \(Q_t^{(n)}\) is the sum of \(n\) independent copies of \(q_t\), the 	classical LLN and CLT apply directly. Moreover,
	\[
	\mathbb{E}[q_t]=g_t'(0),
	\qquad
	\operatorname{Var}(q_t)=g_t''(0),
	\]
	because \(g_t\) is the cumulant generating function of \(q_t\). This proves
	\eqref{eq:lln_sum_iid_queues} and \eqref{eq:clt_sum_iid_queues}.
\end{proof}

\begin{proposition}[LDP for \(Q_t^{(n)}/n\)]
	\label{prop:ldp_sum_iid_queues}
Fix \(t>0\). Let
\(
Q_t^{(n)}:=\sum_{i=1}^{n}q_{i,t},
\)
where \(q_{i,t}\) are independent copies of \(q_t\). Denote the
cumulant generating function of \(q_t\) by
	\[
	g_t(u)
	=
	\log\mathbb{E}\left[e^{uq_t}\right],
	\]
	and let \(\mathcal{D}_{g_t}\) denote its effective domain. Assume that \(g_t\)
	is finite and differentiable on the interior of \(\mathcal{D}_{g_t}\). Then
	\(
	\frac{Q_t^{(n)}}{n}
	\)
	satisfies a large deviation principle with speed \(n\) and  rate function
	\begin{equation}
		I_Q(x)
		=
		\sup_{u\in\mathcal{D}_{g_t}}
		\{ux-g_t(u)\}.
		\label{eq:ldp_sum_iid_queues_rate_function}
	\end{equation}
	Moreover, the  domain of the rate function is 
	$\mathcal{D}_{I_Q}=(0,\infty),$
	and  for \(x\in(0,\infty)\),
	\begin{equation}
		I_Q(x)
		=
		u_xx-g_t(u_x),
		\qquad
		g_t'(u_x)=x.
		\label{eq:ldp_sum_iid_queues_implicit_rate_function}
	\end{equation}
\end{proposition}

\begin{proof}
	As mentioned above independence implies \eqref{eq:BigQ}. By Cramér's theorem \citep{DemboZeitouni2010}, the sequence \(Q_t^{(n)}/n\) satisfies a
	large deviation principle with speed \(n\) and rate function
	\[
	I_Q(x)
	=
	\sup_{u\in\mathcal{D}_{g_t}}
	\{ux-g_t(u)\}.
	\]
		The function \(g_t\) is convex because \(g_t''(u)\geq0\), and
	\(
	\operatorname{Range}(g_t')=(0,\infty).
	\) 
	Note \(g_t\) is finite on all of \(\mathbb{R}\) and convex as proved in Proposition \ref{prop:large_arrival_ldp}.
	Furthermore, for each \(x\in(0,\infty)\), there exists a unique \(u_x\in\mathbb{R}\)
	such that
	\(g_t'(u_x)=x\).
	since \(g_t\) is differentiable and convex, the supremum in the
	Legendre--Fenchel transform is attained at \(u_x\). Therefore,
	\[
	I_Q(x)=u_xx-g_t(u_x),
	\]
	which proves \eqref{eq:ldp_sum_iid_queues_implicit_rate_function}.
\end{proof}

\begin{proposition}[MDP for \(Q_t^{(n)}/n\)]
	\label{prop:mdp_sum_iid_queues}
	Let \((a_n)_{n\geq1}\) be a positive scale satisfying
	\begin{equation}
		\sqrt{n}\ll a_n\ll n,
		\qquad n\to\infty.
		\label{eq:mdp_sum_iid_queues_scale}
	\end{equation}
	Define
	\begin{equation}
		Y_n
		:=
		\frac{Q_t^{(n)}-ng_t'(0)}{a_n}.
		\label{eq:mdp_sum_iid_queues_scaled}
	\end{equation}
	Then \((Y_n)_{n\geq1}\) satisfies a moderate deviation principle with speed
	\(a_n^2/n\) and rate function
	\begin{equation}
		I_Q^{\mathrm{MDP}}(x)
		=
		\frac{x^2}{2g_t''(0)},
		\qquad x\in\mathbb{R}.
		\label{eq:mdp_sum_iid_queues_rate_function}
	\end{equation}
\end{proposition}

\begin{proof}
	For \(\lambda\in\mathbb{R}\), define
	\[
	\Lambda_n(\lambda)
	:=
	\frac{n}{a_n^2}
	\log
	\mathbb{E}\left[
	\exp\left(
	\frac{a_n^2}{n}\lambda Y_n
	\right)
	\right].
	\]
	Since
	\[
	\frac{a_n^2}{n}\lambda Y_n
	=
	\frac{\lambda a_n}{n}Q_t^{(n)}
	-
	\lambda a_n g_t'(0),
	\]
	we have
	\[
	\Lambda_n(\lambda)
	=
	\frac{n^2}{a_n^2}
	g_t\left(\frac{\lambda a_n}{n}\right)
	-
	\frac{\lambda n}{a_n}g_t'(0).
	\]
	Using the Taylor expansion
	\[
	g_t(u)
	=
	g_t(0)+g_t'(0)u+\frac{1}{2}g_t''(0)u^2+o(u^2),
	\qquad u\to0,
	\]
	with \(g_t(0)=0\) and \(u=\lambda a_n/n\), gives
	\[
	\Lambda_n(\lambda)
	\longrightarrow
	\frac{1}{2}g_t''(0)\lambda^2.
	\]
	Since the limiting cumulant generating function is finite and differentiable on
	\(\mathbb{R}\), the Gärtner--Ellis theorem yields the MDP with speed
	\(a_n^2/n\) and rate function
	\[
	\sup_{\lambda\in\mathbb{R}}
	\left\{
	\lambda x-\frac{1}{2}g_t''(0)\lambda^2
	\right\}
	=
	\frac{x^2}{2g_t''(0)}.
	\]
\end{proof}

\section{Numerics}
\label{sec:numerics}

The aim of this section is to illustrate the practical applicability of the
asymptotic results developed in Section~\ref{sec:asymptotics}. In particular,
we numerically evaluate the rate functions associated with the integrated
queue-length process, the queue with infinite arrival intensity, and the sum of
independent and identical queues. These rate functions describe the exponential
decay of probabilities of rare congestion and overload events.

The parameter values used in the numerical experiments are reported in
Table~\ref{tab:asymptotic_parameters}.
\begin{table}[ht]
	\centering
	\caption{Parameter values used in the numerical experiments for the finite-time rate functions.}
	\label{tab:asymptotic_parameters}
	\[
	\begin{array}{ccccccccc}
		\hline
		\text{Parameter} & \kappa & \theta & \sigma & \mu & \lambda_0 & q_0 & t & n \\
		\hline
		\text{Value} & 0.7 & 2 & 0.3 & 1 & 2 & 4 & 10 & 10\\
		\hline
	\end{array}
	\]
\end{table}

\begin{figure}[ht]
	\centering
	\includegraphics[width=0.55\textwidth]{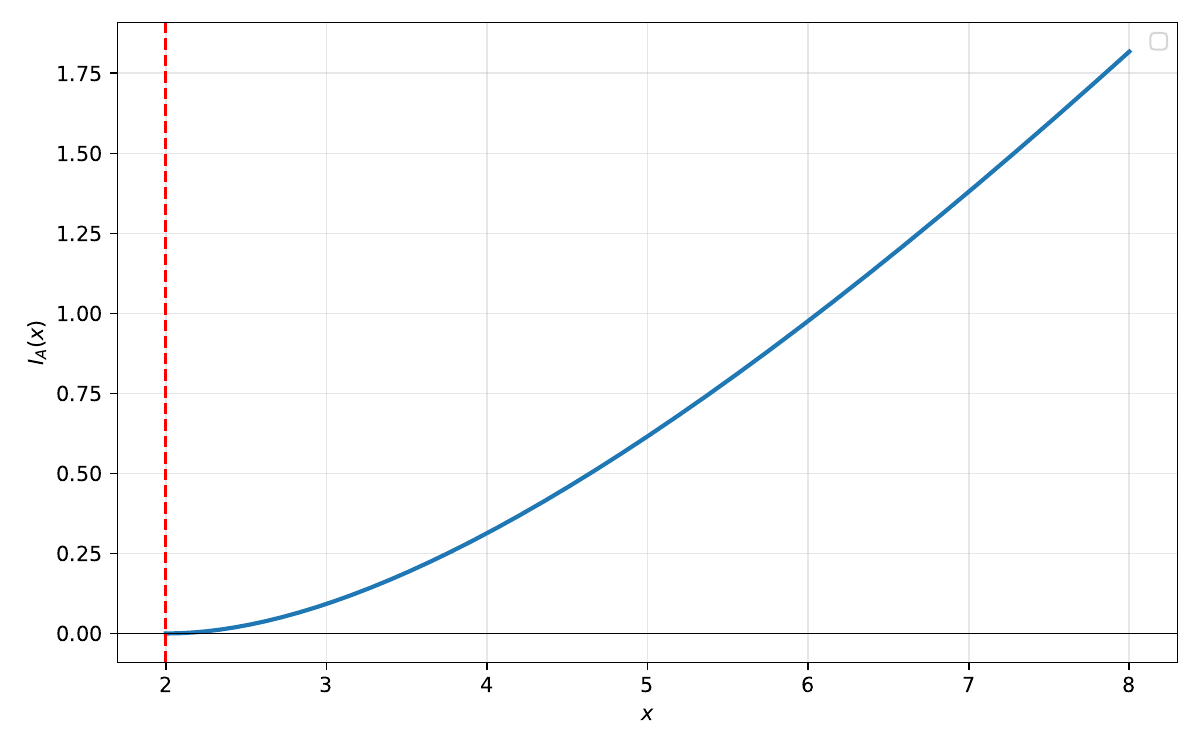}
	\caption{Rate function \(I_A(x)\) for the integrated queue process \(A_t\).}
	\label{fig:rate_fun_integrated_queues}
\end{figure}

We first consider the integrated queue-length process $A_t=\int_0^t q_s\,ds$.
The empirical time average \(A_t/t\) satisfies a large deviation principle with
speed \(t\) and rate function \(I_A(x)\), given in
\eqref{eq:ldp_integrated_queue_rate_function_implicit}. Thus, \(I_A(x)\)
characterises the likelihood of rare deviations of the average queue length
from its typical value. Using the parameter values in
Table~\ref{tab:asymptotic_parameters}, the rate function \(I_A(x)\) is computed
numerically. Figure~\ref{fig:rate_fun_integrated_queues} displays the resulting
rate function in blue. The vertical dashed line marks the LLN point
$x=\frac{\theta}{\mu},$
where the rate function attains its minimum value \(0\).

\begin{figure}[ht]
	\centering
	\subfloat[]
	{\includegraphics[width=0.5\textwidth,height=5.5cm]{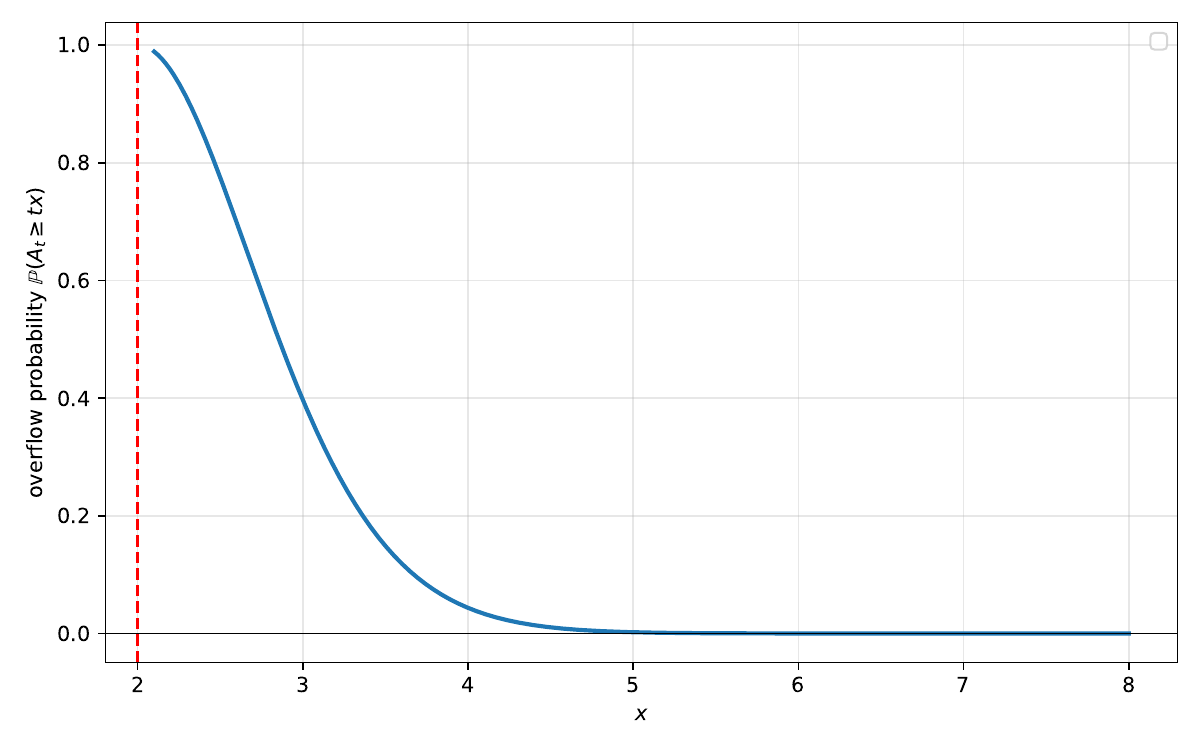}}
	\subfloat[]
	{\includegraphics[width=0.5\textwidth,height=5.5cm]{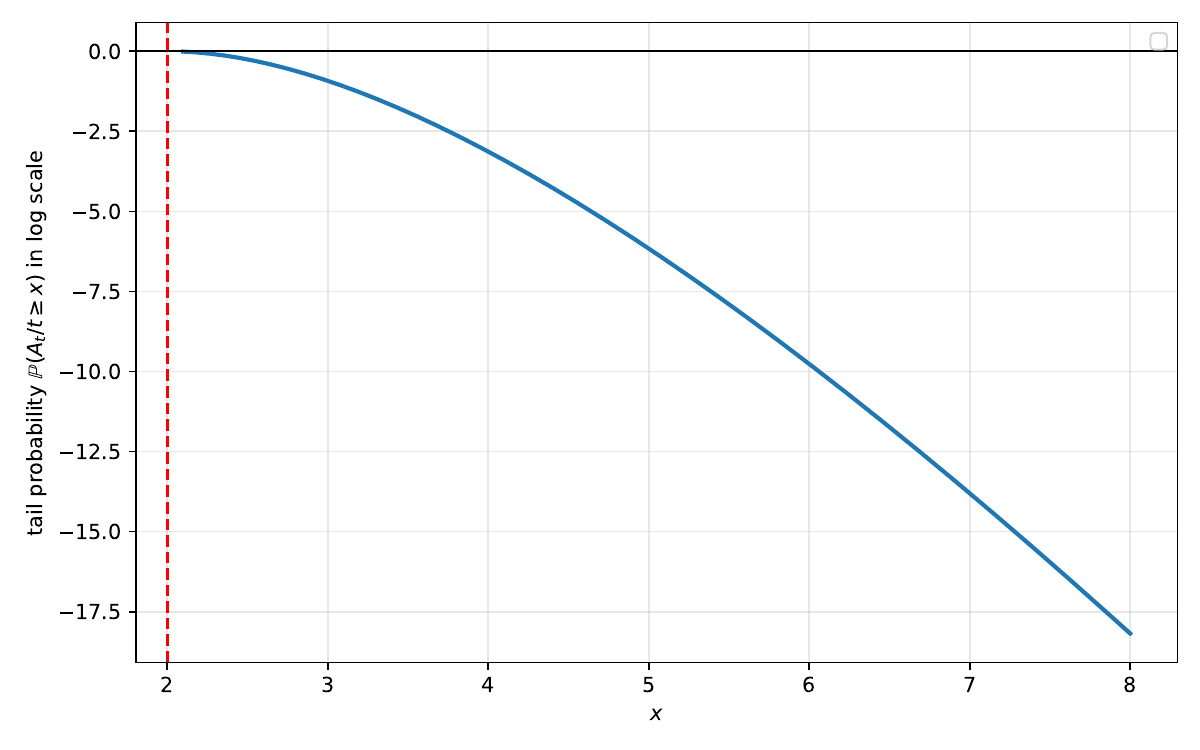}}
	\caption{(a) Overflow probability approximation based on \(I_A(x)\); (b)
		tail probability approximation in logarithmic scale based on \(I_A(x)\).}
	\label{fig:prob_integrated_queues}
\end{figure}

For the integrated queue-length process \(A_t\) at the fixed time \(t=10\), and
for \(x\) larger than the typical level \(x=\theta/\mu\), the overflow
probability associated with congestion and overload events is approximated by
\begin{equation}
	\mathbb{P}\left(A_t \geq tx\right)
	\sim
	e^{-t I_A(x)}.
	\label{eq:overflow_integrated_queues}
\end{equation}
This approximation is illustrated in
Figure~\ref{fig:prob_integrated_queues}(a). Similarly, the upper-tail
asymptotics are governed by the same rate function:
\begin{equation}
	\log\mathbb{P}\left(\frac{A_t}{t}>x\right)
	\sim
	-tI_A(x).
	\label{eq:tail_integrated_queues}
\end{equation}
Figure~\ref{fig:prob_integrated_queues}(b) shows the corresponding 
decay of the tail probabilities on a logarithmic scale.

\begin{figure}[ht]
	\centering
	\includegraphics[width=0.55\textwidth]{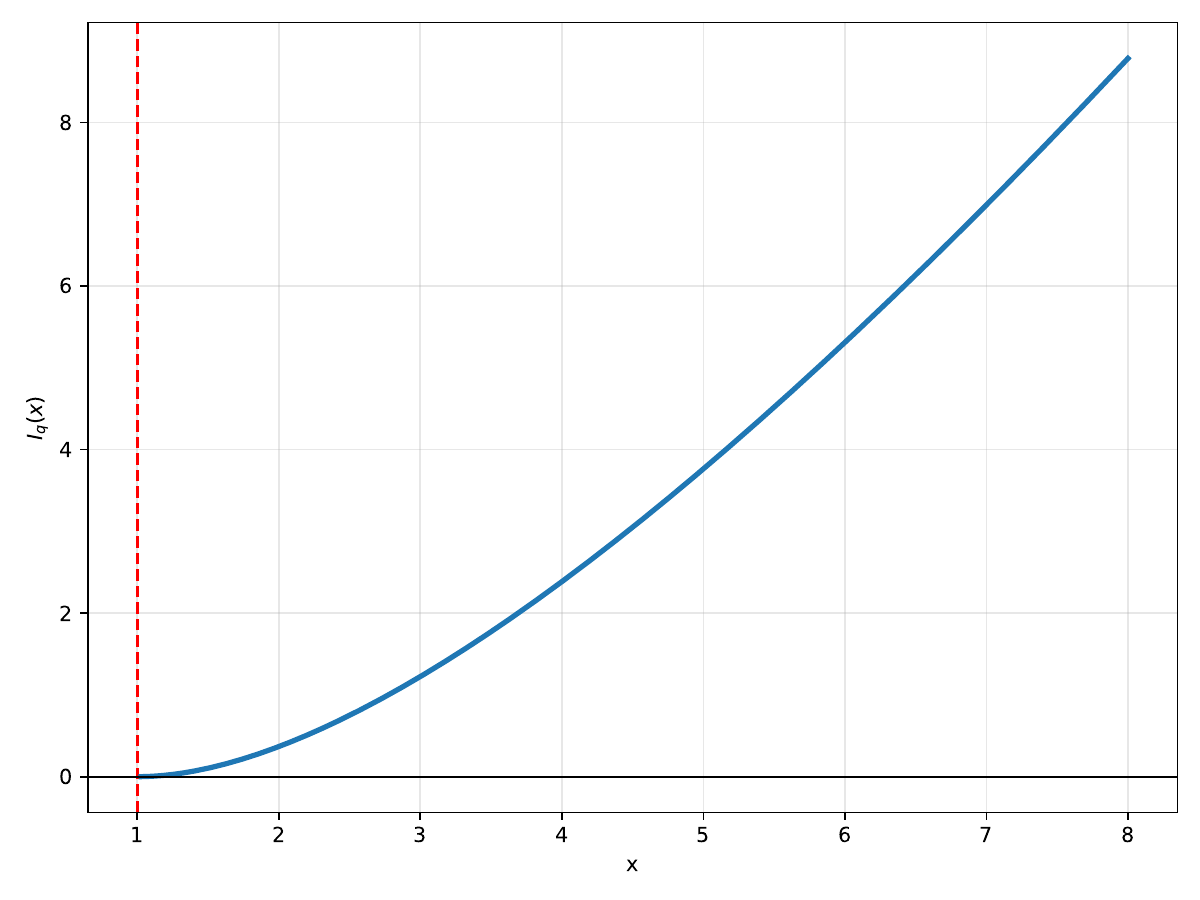}
	\caption{Rate function \(I_q(x)\) for the queue process with infinite arrival intensity
		\(q_t^{(n)}\).}
	\label{fig:rate_fun_infinite_queues}
\end{figure}

We next consider the queue process with infinite arrival intensity,
\(q_t^{(n)}\). The scaled queue length \(q_t^{(n)}/n\) satisfies a large
deviation principle with speed \(n\) and rate function \(I_q(x)\), given in
\eqref{eq:large_arrival_ldp_rate_function_implicit}. Therefore, \(I_q(x)\)
characterises rare deviations of the scaled queue length from its limiting
mean. Using the parameter values in Table~\ref{tab:asymptotic_parameters}, the
rate function \(I_q(x)\) is computed numerically. Figure~\ref{fig:rate_fun_infinite_queues}
displays the resulting rate function in blue. The vertical dashed line marks
the LLN point
\[
x=
\frac{\kappa\mu}{\kappa-\mu}
\left[
\frac{1-e^{-\mu t}}{\mu}
-
\frac{1-e^{-\kappa t}}{\kappa}
\right],
\]
where the rate function attains its minimum value \(0\).

\begin{figure}[ht]
	\centering
	\subfloat[]
	{\includegraphics[width=0.5\textwidth,height=5.5cm]{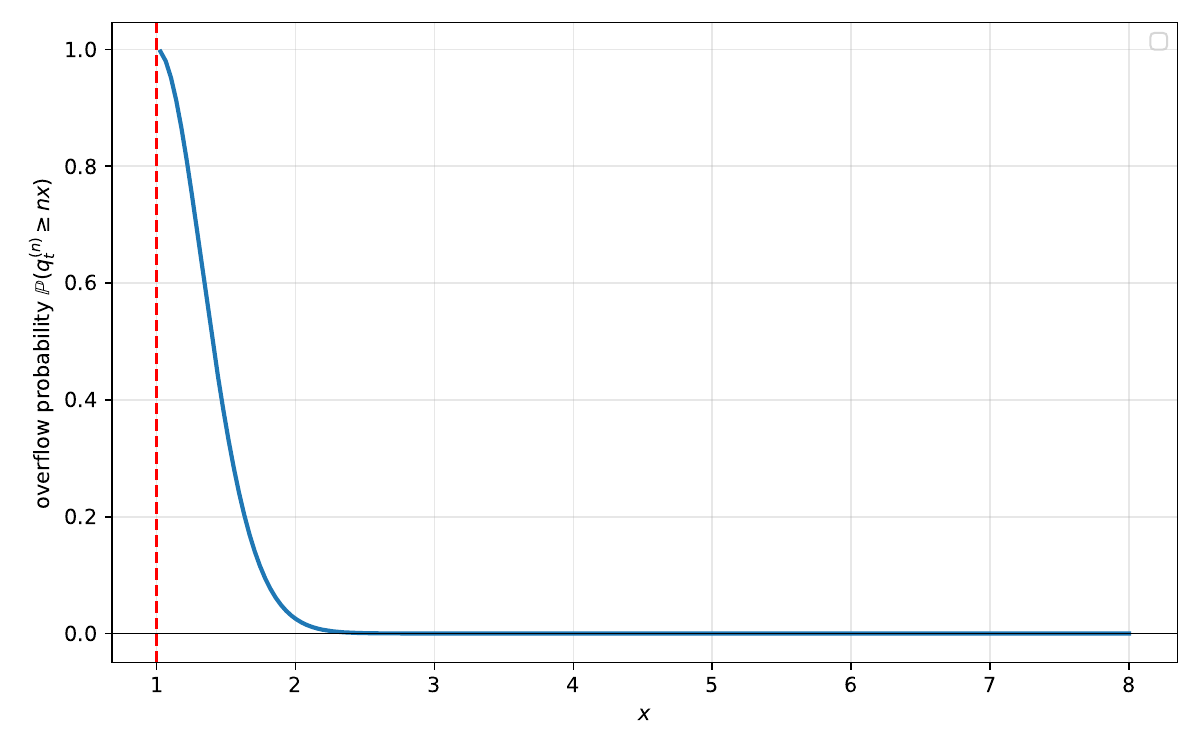}}
	\subfloat[]
	{\includegraphics[width=0.5\textwidth,height=5.5cm]{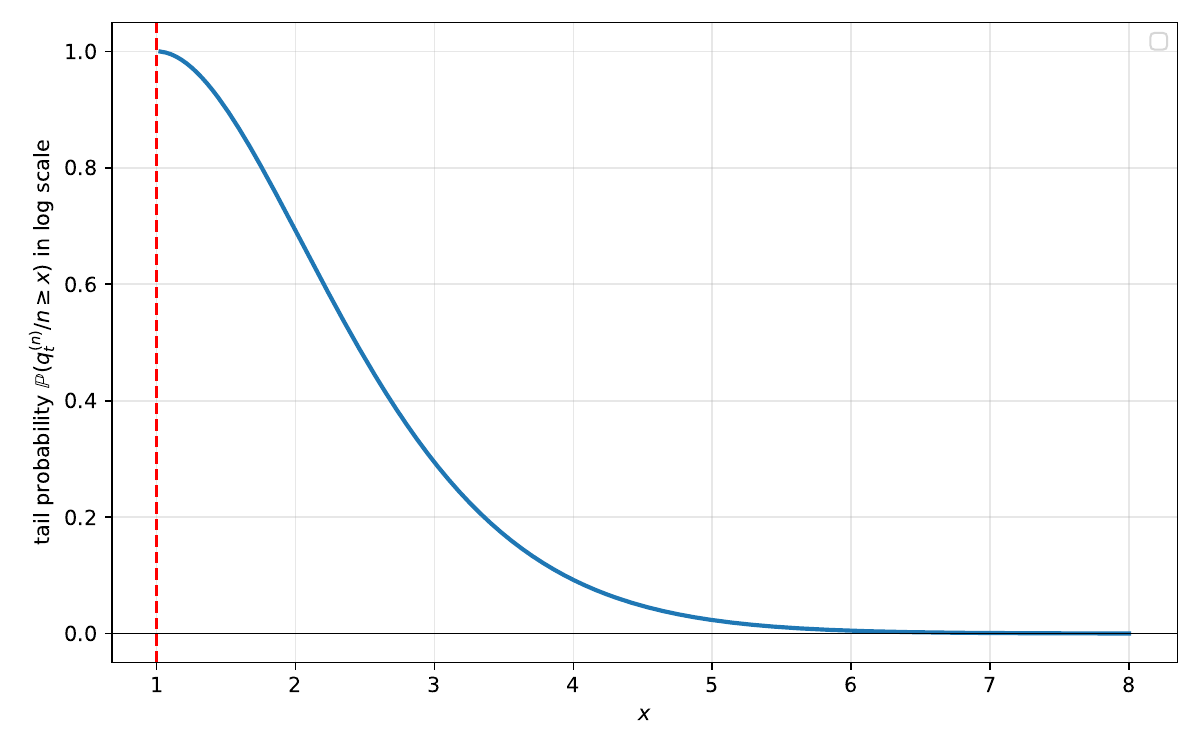}}
	\caption{(a) Overflow probability approximation based on \(I_q(x)\); (b)
		tail probability approximation in logarithmic scale based on \(I_q(x)\).}
	\label{fig:prob_infinite_queues}
\end{figure}

For the queue process with infinite arrival intensity \(q_t^{(n)}\), with fixed
\(n=10\), and for \(x\) larger than the typical level
\[
x=
\frac{\kappa\mu}{\kappa-\mu}
\left[
\frac{1-e^{-\mu t}}{\mu}
-
\frac{1-e^{-\kappa t}}{\kappa}
\right],
\]
the overflow probability associated with congestion and overload events is
approximated by
\begin{equation}
	\mathbb{P}\left(q_t^{(n)} \geq nx\right)
	\sim
	e^{-n I_q(x)}.
	\label{eq:overflow_infinite_queues}
\end{equation}
This approximation is illustrated in Figure~\ref{fig:prob_infinite_queues}(a).
Similarly, the upper-tail asymptotics are controlled by the same rate function:
\begin{equation}
	\log\mathbb{P}\left(\frac{q_t^{(n)}}{n}>x\right)
	\sim
	-nI_q(x).
	\label{eq:tail_infinite_queues}
\end{equation}
Figure~\ref{fig:prob_infinite_queues}(b) shows the decay of the
tail probabilities on a logarithmic scale.

\begin{figure}[ht]
	\centering
	\includegraphics[width=0.55\textwidth]{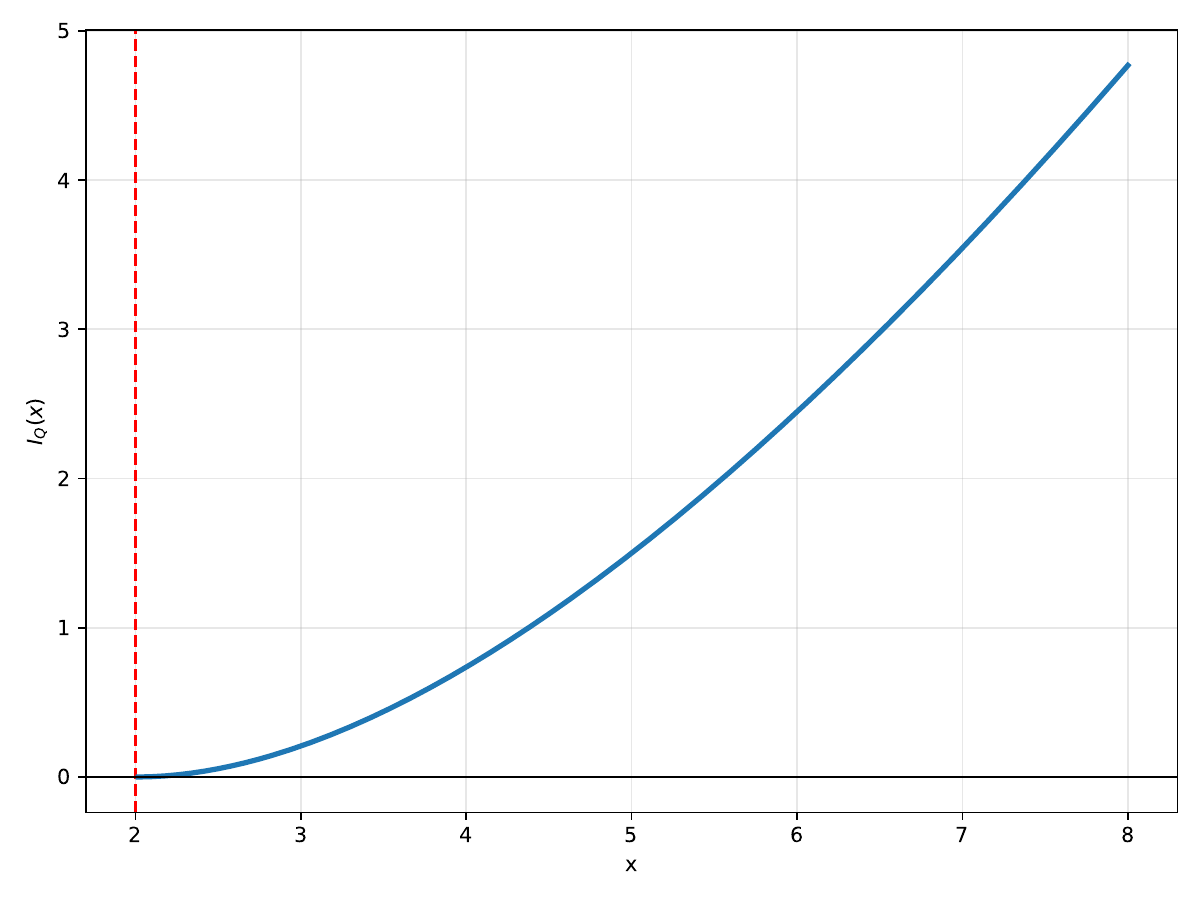}
	\caption{Rate function \(I_Q(x)\) for the sum process \(Q_t^{(n)}\) of \(n\)
		identical and independent queues.}
	\label{fig:rate_fun_sum_queues}
\end{figure}

Finally, we consider the sum process \(Q_t^{(n)}\) of \(n\) identical and
independent queues $q_t$. The empirical average \(Q_t^{(n)}/n\) satisfies a large
deviation principle with speed \(n\) and rate function \(I_Q(x)\), given in
\eqref{eq:ldp_sum_iid_queues_rate_function}. Hence, \(I_Q(x)\) characterises
rare deviations of the average queue length from its typical level. Using the
parameter values in Table~\ref{tab:asymptotic_parameters}, the rate function
\(I_Q(x)\) is computed numerically. Figure~\ref{fig:rate_fun_sum_queues}
displays the resulting rate function in blue. The vertical dashed line marks
the LLN point
\(
x=\frac{\theta}{\mu},
\)
where the rate function attains its minimum value \(0\).

\begin{figure}[ht]
	\centering
	\subfloat[]
	{\includegraphics[width=0.5\textwidth,height=5.5cm]{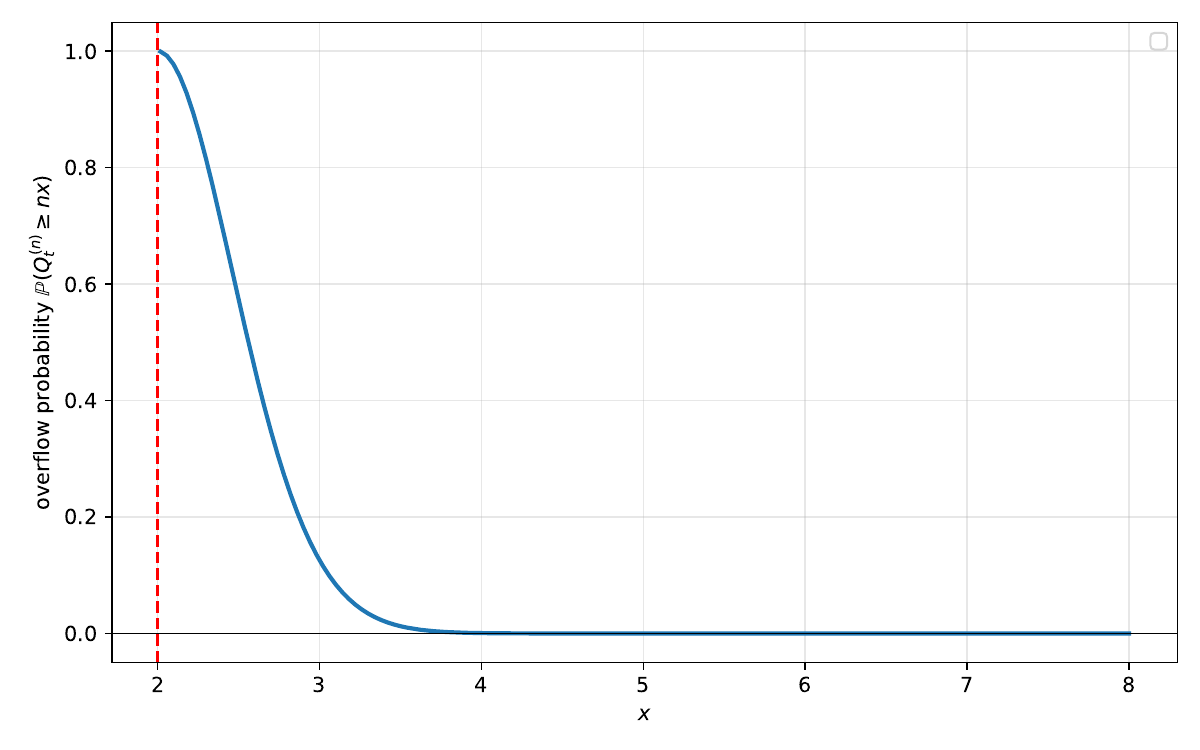}}
	\subfloat[]
	{\includegraphics[width=0.5\textwidth,height=5.5cm]{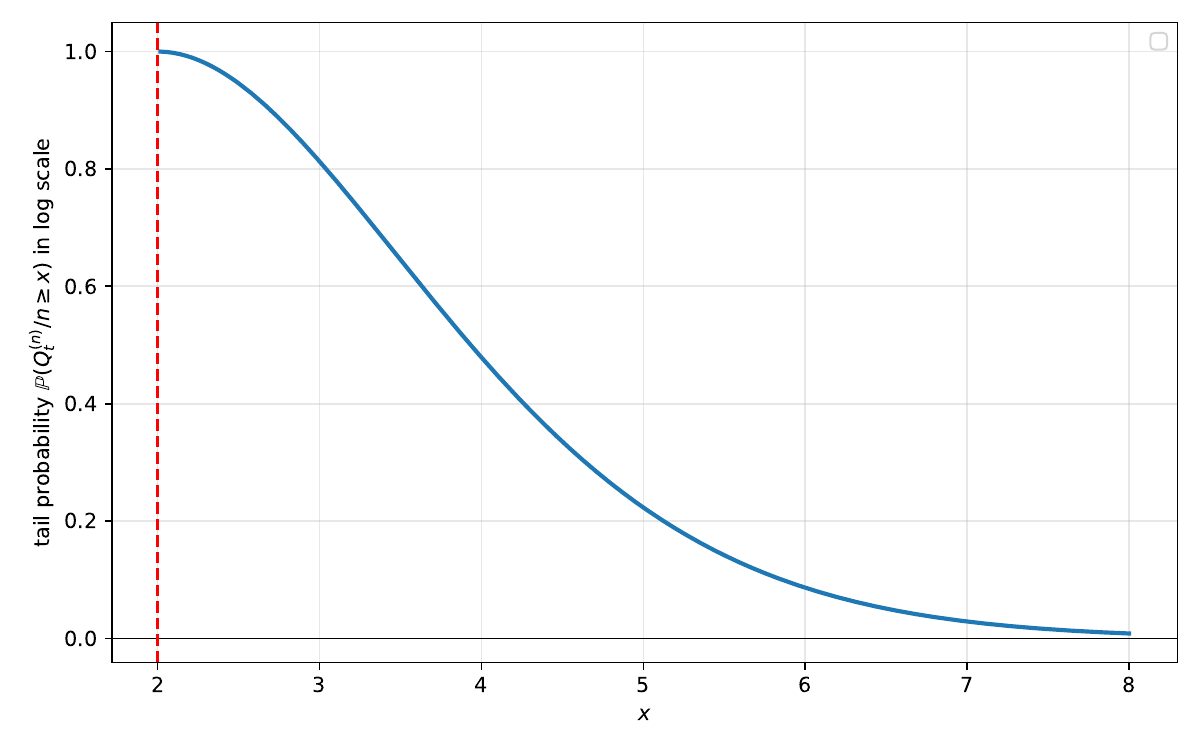}}
	\caption{(a) Overflow probability approximation based on \(I_Q(x)\); (b)
		tail probability approximation in logarithmic scale based on \(I_Q(x)\).}
	\label{fig:prob_sum_queues}
\end{figure}

For the sum process \(Q_t^{(n)}\), with fixed \(n=10\), and for \(x\) larger
than the typical level \(x=\theta/\mu\), the overflow probability associated
with congestion and overload events is approximated by
\begin{equation}
	\mathbb{P}\left(Q_t^{(n)} \geq nx\right)
	\sim
	e^{-n I_Q(x)}.
	\label{eq:overflow_sum_queues}
\end{equation}
Similarly, the upper-tail asymptotics are governed by the same rate function:
\begin{equation}
	\log\mathbb{P}\left(\frac{Q_t^{(n)}}{n}>x\right)
	\sim
	-nI_Q(x).
	\label{eq:tail_sum_queues}
\end{equation}
Figure~\ref{fig:prob_sum_queues}(a) illustrates the overflow probability
approximation, while Figure~\ref{fig:prob_sum_queues}(b) shows the exponential
decay of the tail probabilities on a logarithmic scale.

\section{Conclusion}\label{sec:Conclusion}

This paper proposes a model for infinite-server queues driven by a Cox--BESQ process that admits explicit formulas for the joint moment-generating function (MGF) of the arrival intensity, the number of customers in the system, and the cumulative number of customers in the system over a given time interval. These results are made possible by the affine and tractable structure of the model. Using these MGFs, we derive several asymptotic results, including the LLN, CLT, LDP, and MDP, for various quantities associated with the Cox--BESQ infinite-server queue. More precisely, we consider the cumulative number of customers in the system over a time interval, an infinite-server queue with a diverging arrival intensity, and sums of independent infinite-server queues. For each of these three important cases, we establish all four asymptotic results and illustrate their practical applicability through numerical examples.

Several extensions are possible. One direction is to consider more general service-time distributions beyond the assumptions adopted in the present analysis. Another is to consider more general arrival-intensity dynamics. For example, incorporating jumps while preserving the affine and tractable structure of the model would be of considerable interest. It would also be valuable to investigate statistical estimation procedures for the model parameters and to calibrate the model to empirical arrival data. Finally, empirical applications to staffing and resource-allocation problems would further enhance the practical relevance of the proposed framework.

\clearpage

\footnotesize
\bibliographystyle{abbrvnat}
\bibliography{Biblio}
\end{document}